\documentclass[12pt,a4paper]{article}

\usepackage[left=1cm,right=1cm,top=2cm,bottom=2cm]{geometry}
\usepackage[dvips]{graphics}
\usepackage[dvips]{epsfig}
\usepackage{color}
\usepackage{hyperref}
\usepackage{tikz,authblk}

\newtheorem{lemma}{Lemma}
\newtheorem{theorem}[lemma]{Theorem}
\newtheorem{corollary}[lemma]{Corollary}
\newtheorem{proposition}[lemma]{Proposition}
\newtheorem{definition}{Definition}
\newtheorem{remark}{Remark}
\newtheorem{example}{Example}

\newcommand{\dimo}{\noindent \emph{Proof. }}
\newcommand{\qed}{\\ \rightline{$\Box$ \ \ \ \ \ \ \ \ \ \ \ \ \ \ \ }\\}
\newcommand{\e}{\varepsilon}
\newcommand{\G}{\Gamma}
\newcommand{\g}{\gamma}

\usepackage{latexsym}
\usepackage{amsfonts}
\usepackage{amssymb}
\usepackage{amsmath}
\usepackage[english]{babel}

\begin{document}

\title{Kirby diagrams, trisections and gems of PL 4-manifolds: \\ 
relationships, results and open problems}

 \renewcommand{\Authfont}{\scshape\small}
 \renewcommand{\Affilfont}{\itshape\small}
 \renewcommand{\Authand}{ and }

\author[1] {Maria Rita Casali}
\author[2] {Paola Cristofori}

\affil[1,2] {Department of Physics, Informatics and Mathematics, University of Modena and Reggio Emilia\par  Via Campi 213 B, I-41125 Modena (Italy), casali@unimore.it,\ paola.cristofori@unimore.it}

\maketitle

 \abstract{We review the main achievements regarding the interactions between gem theory (which makes use of edge-colored graphs to represent PL-manifolds of arbitrary dimension) and both the classical representation of PL 4-manifolds via Kirby diagrams and the more recent one via trisections. 
Original results also appear (in particular, about gems representing closed 4-manifolds which need 3-handles in their handle decomposition, as well as about trisection diagrams), together with open problems and further possible applications to the study of compact PL 4-manifolds. }
\endabstract


  

\section{Introduction}

Gem theory is a well-established representation theory for PL-manifolds of arbitrary dimension, whose combinatorial tool are edge-colored graphs, dual to colored triangulations (see, for example, \cite{generalized-genus} and references therein). 
In dimension 4, it has been proved to have interesting interactions both with the classical representation via Kirby diagrams (or, equivalently, handle decompositions) and with the more recent one via trisections. In particular, algorithmic procedures have been obtained, to pass from one representation to the other (see \cite{Casali-Cristofori Kirby-diagrams}, \cite{Casali-Cristofori gem-induced} and \cite{Casali-Cristofori trisection bis}), allowing to yield estimations of the involved invariants and to produce examples of triangulations of exotic 4-manifolds in order to investigate their PL structures.

The present paper reviews the main achievements on the subject, together with open problems and further possible applications to the study of compact PL 4-manifolds. In order to help the reader gain a deep understanding of the above interactions, without forcing him to continually refer to previous papers, sketches of the main proofs are provided, so as to make the text essentially self-contained. 

Original results also appear, which improve the state-of-the-art and open to further developments. In particular, some typical situations are identified, where it is possible to obtain a gem of a closed 4-manifold described by a handle decomposition, starting from a gem of the compact 4-manifold consisting only of the handles up to index two: see Proposition \ref{chiudere} and related examples.       
Moreover, in Section \ref{ss.from_Kirby_to_trisection_diagrams}, we explain how to obtain trisection diagrams directly from gems admitting gem-induced trisections and/or from Kirby diagrams, via gem theory. 
\medskip

The paper is organized as follows. 

Section 2 is devoted to a brief review of the basic notions about the three involved representation methods for PL 4-manifolds, i.e. Kirby diagrams, trisections and colored graphs.   
Section 3 contains the main results and constructions about the relationships among the three theories: in particular, from Kirby diagrams to gems (Section \ref{ss.from_Kirby_to_gems}, which contains also the new result related to the addition of 3-handles), from gems to trisections (Section \ref{ss.from_gems_to_trisections}) and from Kirby diagrams to trisections via gems (Section \ref{ss.from_Kirby_to_trisections}).  Within the latter, the original part regarding trisection diagrams is presented, showing how trisection diagrams can be easily obtained from gems of 4-manifolds admitting gem-induced trisections (Proposition \ref{trisection_diagrams_gem-induced}), and in particular from gems associated to Kirby diagrams (Proposition \ref{trisection_diagrams_framed-links}); by a suitable extension of the notion of trisection diagram, similar results are also obtained for simply-connected compact 4-manifolds with connected boundary (Proposition \ref{trisection_diagrams_simply-connected_bounded}), in particular starting from their Kirby diagrams (Proposition \ref{trisection_diagrams_gem-induced_boundary}). All these results closely resemble the already known ones in dimension 3, where gems of 3-manifolds, both in the closed and boundary case, directly yield 
Heegaard diagrams.   

Finally, Section 4 completes the paper with a list of open problems, suggesting trends for future investigation. 


\section{Preliminaries about the involved representations}
In this Section we briefly review the basic notions of the three representation methods for PL 4-manifolds involved in the paper, without any claim to exhaustiveness. Throughout the paper all manifolds and maps, when not specified, are understood to be PL. For generalities about 4-dimensional topology of piecewise linear manifolds, we refer to \cite{[M]} and \cite{GS}. 
\medskip

In the following,  ${\mathbb Y}_m$ will denote the orientable 4-dimensional handlebody of genus $m$, i.e. $\natural_m (\mathbb S^1 \times \mathbb D^3)$, while  
${\mathbb Y}_m^{(\sim)}$ (resp. $\mathbb S^1\otimes\mathbb S^{n-1}$) will denote a  genus $g$ 4-dimensional handlebody (resp.  an $\mathbb S^{n-1}$-bundle over $\mathbb S^1$), orientable or non-orientable according to the context. 

Obviously ${\partial{\mathbb Y}}_m^{(\sim)} = \#_m(\mathbb S^1 \otimes \mathbb S^2)$ and, whenever orientability is assumed, the $n$-sphere (resp. the $4$-disk) will be included in the notation $\#_m(\mathbb S^1 \otimes \mathbb S^{n-1}) = \#_m(\mathbb S^1 \times \mathbb S^{n-1})$  (resp. ${\mathbb Y}_m^{(\sim)}={\mathbb Y}_m$) as the case $m=0.$


\subsection{Kirby diagrams of orientable $4$-manifolds}\label{ss:Kirby diagrams} 
\ \ 

It is well-known that, if a compact orientable 4-manifold $M$ admits a handle decomposition with no 3-handles, then $M$ can be represented by a {\it Kirby diagram}, i.e. a link in the 3-sphere equipped with some further information, which encodes the attachments of 1- and 2-handles.

In fact, since the existence of exactly one $0$-handle can always be assumed, then 1- and 2-handles can be visualized, by positioning on the boundary of the $0$-handle, as follows:
\begin{itemize}
    \item the attachment of each 2-handle is represented by a {\it knot} (i.e. the core of its attaching solid torus) equipped with an integer specifying the {\it framing} of the handle (i.e. the choice of a normal vector field on the knot); 
    \item the attachment of a 1-handle would be specified by its attaching region formed by two 3-balls. However, this representation involves some technical issues, which lead to often prefer the so-called {\it Akbulut's notation}: this is based on the observation that the attachment of a 1-handle can be also thought of as the drilling of its complementary 2-handle. Therefore, each 1-handle can be visualized by drawing the boundary of the cocore of its complementary 2-handle, which is an unknot with no framing, and it is marked by a dot to distinguish it from the components that come from ``real" 2-handles  (\cite{GS}).
\end{itemize}

The resulting object is therefore a link $L$ in the 3-sphere having some {\it framed} components equipped with integers and some {\it dotted} ones. Furthermore, it is always possible to arrange the handles in such a way that the dotted components, if any, are pairwise unlinked.  

Obviously $M$ can be re-constructed from $\mathbb D^4$ by adding 1-handles according to the dotted and 2-handles according to the framed components of $L$; furthermore, $\partial M$ is the closed 3-manifold obtained by Dehn surgery on the framed link obtained by replacing each dotted component of $L$ by a 0-framed one (\cite{GS}  and \cite{[M]}). 

In particular, if $\partial M\cong\#_r(\mathbb S^1\times\mathbb S^2)$ ($r\geq 0$), a famous result by Laudenbach and Poenaru (\cite{Laudenbach-Poenaru}) ensures that there is a unique way to attach $r$ 3-handles and a 4-handle to obtain a closed manifold $\bar M$, which is then uniquely associated with $M$.
As a consequence, in this case, any Kirby diagram representing $M$ can be thought of as representing $\bar M = M \cup \mathbb Y_r$  as well.

In \cite{Montesinos}, Kirby diagrams representing closed orientable $4$-manifolds (i.e. representing compact $4$-manifolds with boundary homeomorphic to $\#_r(\mathbb S^1\times\mathbb S^2),$ with $r\geq 0$) are also called {\it (4-dimensional) Heegaard diagrams}.  


\subsection{Trisections of closed $4$-manifolds}  \label{ss: trisections} 

\ \ 

By generalizing the classical idea of Heegaard splitting in dimension 3, in 2016 Gay and Kirby (\cite{Gay-Kirby}) introduced the notion of {\it trisection} of a smooth closed oriented $4$-manifold; then, the notion was extended to non-orientable 4-manifolds in \cite{Spreer-Tillmann(Exp)} and \cite{Miller-Naylor}.

\begin{theorem}  \ {\rm (\cite{Gay-Kirby}, \cite{Spreer-Tillmann(Exp)}, \cite{Miller-Naylor})} \  \label{th: trisection}
Each smooth closed orientable (resp. non-orientable) $4$-manifold $M$ admits a decomposition $M = H_0 \cup H_1 \cup H_2$, such that:  
 \begin{itemize}
 \item [(a)]  $H_{0}, H_{1}, H_{2}$ are $4$-dimensional orientable (resp. non-orientable) handlebodies with pairwise disjoint interiors; 
 \item [(b)]  $H_{01}=H_{0}\cap H_{1},$  $H_{02}=H_{0}\cap H_{2}$ and $H_{12}=H_1\cap H_2$ are $3$-dimensional orientable (resp. non-orientable) handlebodies with the same genus; 
 \item [(c)]  $\Sigma = H_{0}\cap H_{1}\cap H_{2}$  is a closed connected orientable (resp. non-orientable) surface. 
 \end{itemize}
\end{theorem} 

\begin{definition} \ {\rm (\cite{Gay-Kirby}, \cite{Spreer-Tillmann(Exp)}, \cite{Miller-Naylor})} \  \label{def: trisection}
{\em A triple $(H_0, H_1, H_2)$ that satisfies the requirements of Theorem \ref{th: trisection} is said to be a {\it trisection} of $M$; the {\it genus} of such a trisection is defined as the common genus of the 3-dimensional handlebodies $H_{01}, H_{02}, H_{12}$, while the surface $\Sigma$ is called the {\it central surface} of the trisection.\footnote{Note that in the orientable (resp. non-orientable) case the genus of the central surface equals (resp. is the double of) the genus of the trisection.  Actually, many slightly different versions of the notions regarding trisections exist, also in the cited papers, for example requiring that all $4$-dimensional handlebodies have the same genus (giving rise to the so-called {\it balanced trisections}), or considering the genus of the central surface as the genus of the trisection, also in the non-orientable case. However, no loss of generality occurs, and it is not difficult to translate one version into the others.} 

The {\it trisection genus} $g_T(M)$ of a smooth closed $4$-manifold $M$ is defined as the minimum genus of a trisection of $M$.  }
\end{definition}  

\smallskip 

It is not difficult to check that, if $(H_0, H_1, H_2)$ is a genus $g$ trisection of $M$, then $(\Sigma, H_{ij}, H_{ik})$ is a genus $g$ Heegaard splitting of $\partial H_i$, for each permutation $(i,j,k)$ of $\{0,1,2\}$.  
Moreover,  if $k_i$ ($i \in \{0,1,2\}$) denotes the genus of the $4$-dimensional ``piece" $H_i$, the equality $\chi(M)= 2+g-k_0-k_1-k_2$ holds, together with the inequalities $g \ge \beta_1(M; \mathbb Z_2) + \beta_2(M; \mathbb Z_2)$ and $\beta_1(M; \mathbb Z_2) \le k_i$, $\forall i \in \{0,1,2\}$. 

As a consequence, only simply-connected $4$-manifolds can admit trisections where one the $4$-dimensional ``pieces" is a $4$-disk. Viceversa, it is an open problem whether each closed simply-connected $4$-manifold admits a trisection with this property, or not: see, for example, \cite{Meier}, \cite{Lambert-Cole-Meier} and \cite{Meier-Schirmer-Zupan}. 

\bigskip
In virtue of the already cited theorem by Laudenbach and Poenaru (\cite{Laudenbach-Poenaru}), together with its non-orientable version proved in \cite{Miller-Naylor}, the attachments of the 3-dimensional handlebodies on the central surface are sufficient to identify the trisection, and hence the represented closed $4$-manifold: 
 
\begin{definition} \label{def. trisection-diagram}{\em A $(g; k_0; k_1; k_2)$-{\it trisection diagram} is a 4-tuple $(\Sigma; \alpha, \beta,\gamma)$ such that the triples $(\Sigma; \alpha, \beta)$, $(\Sigma; \alpha, \gamma)$ and $(\Sigma; \beta,\gamma)$ are genus $g$ Heegaard diagrams for $\#_{k_0}(\mathbb S^1 \otimes \mathbb S^2)$, $\#_{k_1}(\mathbb S^1 \otimes \mathbb S^2)$ and $\#_{k_2}(\mathbb S^1 \otimes \mathbb S^2)$ respectively.\footnote{Obviously, all $3$- and $4$-dimensional manifolds involved in Definition \ref{def. trisection-diagram}, as well as in the procedure to reconstruct the represented closed $4$-manifold, are orientable or not according to $\Sigma$.} }
\end{definition}  

In fact - as explained in \cite{Miller-Naylor} both in the orientable and non-orientable setting - the $4$-manifold $M$ is uniquely obtained, up to diffeomorphism, in the following way: 
 \begin{itemize}
\item[-]  take the product $\Sigma \times \mathbb D^2$; 
  \item[-] take three genus $g$ 3-dimensional handlebodies $V_\alpha$, $V_\beta$ and $V_\gamma$, identified by the systems of curves $\alpha,$ $\beta$ and $\gamma$ respectively; 
  \item[-]  attach $V_\alpha \times I$ (resp. $V_\beta \times I)$ (resp. $V_\gamma \times I)$ to $\partial (\Sigma \times \mathbb D^2) = \Sigma \times \mathbb S^1$ along $\Sigma \times [-\epsilon, \epsilon]$ (resp. $\Sigma \times [\frac 2 3 \pi -\epsilon, \frac 2 3 \pi + \epsilon]$)  (resp. $\Sigma \times [\frac 4 3 \pi -\epsilon, \frac 4 3 \pi + \epsilon]$);  
  \item[-]  fill the three boundary components (isomorphic to $\#_{k_i}(\mathbb S^1 \otimes \mathbb S^2)$, for $i=0,1,2$) of the resulting compact $4$-manifold with three $4$-dimensional handlebodies (of genus $k_i$, for $i=0,1,2$).       
\end{itemize}

\medskip 

\subsection{Gems} \label{preliminaries.gems}   


{\it Crystallization theory}, or {\it gem theory}, is a representation tool for piecewise linear (PL) compact manifolds, without assumptions about dimension, connectedness, orientability or boundary properties. 
For a detailed review of the foundational notions and results, see the ``classical'' survey paper \cite{Ferri-Gagliardi-Grasselli}, or the book 
\cite{Lins-book} regarding dimension $3$; for updated developments, especially in dimension $4$, see \cite{generalized-genus} and \cite{Casali-Cristofori-Gagliardi RIMUT 2020}, together with their references. 
Here, we summarize only the basic elements that are necessary for the topic we are interested in. 


\begin{definition} \label{$n+1$-colored graph}{\em An $(n+1)${\emph{-colored graph}}  ($n \ge 2$) is a pair $(\G,\g)$, where $\G=(V(\G), E(\G))$ is a multigraph (i.e. multiple edges are allowed, while loops are forbidden)  which is regular of degree  $n+1$, and $\g$ is an {\it edge-coloration}, that is a map  $\g: E(\G) \rightarrow \Delta_n=\{0,\ldots, n\}$ which is injective on adjacent edges.

For sake of concision, the coloration $\gamma$ is often understood, and the colored graph is simply denoted by $\G$.} \end{definition}


For every  $\{c_1, \dots, c_h\} \subseteq\Delta_n$ let $\G_{\{c_1, \dots, c_h\}}$  be the subgraph obtained from $\G$  by deleting all the edges that are not colored by the elements of $\{c_1, \dots, c_h\}$. 
In this setting, the complementary set of $\{c\}$ (resp. $\{c_1,\dots,c_h\}$)  in $\Delta_n$ will be denoted by $\hat c$ (resp. $\hat c_1\cdots\hat c_h$). 
The connected components of $\G_{\{c_1, \dots, c_h\}}$ are called {\it $\{c_1, \dots, c_h\}$-residues} or {\it $h$-residues} of $\G$; their number is indicated by $g_{\{c_1, \dots, c_h\}}$ (or, for short, by $g_{c_1,c_2}$, $g_{c_1,c_2,c_3}$ and $g_{\hat c}$ if $h=2,$ $h=3$ and $h = n$ respectively). 

 \medskip 

Within gem theory, an $(n+1)$-colored graph $\G$ is thought of as the combinatorial visualization of an $n$-dimensional pseudocomplex $K(\G)$, which is associated to $\Gamma$ in the following way:
\begin{itemize}
\item take an $n$-simplex for each vertex of $\G$ and label its vertices by the elements of $\Delta_n$;
\item if two vertices of $\G$ are $c$-adjacent ($c\in\Delta_n$), glue the corresponding $n$-simplices  along their $(n-1)$-dimensional faces opposite to the $c$-labelled vertices, so that equally labelled vertices are identified.
\end{itemize}

\smallskip
In general, $\vert K(\G)\vert$ is an {\it n-pseudomanifold} and $\G$ is said to {\it represent} it. 

\medskip

By construction, $K(\G)$ is endowed with a vertex-labeling by $\Delta_n$ that is injective on each simplex. Moreover, $\G$ turns out to be the 1-skeleton of the dual complex of $K(\G)$.
The duality establishes a bijection between the $\{c_1, \dots, c_h\}$-residues of  $\G$  
and the $(n-h)$-simplices of $K(\G)$ whose vertices are labelled by $\Delta_n - \{c_1, \dots, c_h\}$. 
In particular, given an $(n+1)$-colored graph $\G$, each connected component of $\G_{\hat c}$ ($c\in\Delta_n$) is an $n$-colored graph representing the link of a $c$-labelled vertex of $K(\G)$ in the first barycentric subdivision of $K(\G).$ 

\noindent 

As a consequence, the following characterizations hold: 

\begin{itemize}
    \item $\vert K(\G)\vert$ is a closed $n$-manifold iff, for each color $c\in\Delta_n$, all $\hat c$-residues of $\G$ represent the $(n-1)$-sphere; 
    \item  $\vert K(\G)\vert$ is a singular\footnote{A {\it singular (PL) $n$-manifold} is a closed connected $n$-dimensional polyhedron admitting a simplicial triangulation where the links of vertices are closed connected $(n-1)$-manifolds, while the links of $h$-simplices, for $h \ge 1$, are $(n-h-1)$-spheres. Vertices whose links are not $(n-1)$-spheres are called {\it singular}. The notion extends also to polyhedra associated to colored graphs, provided that links are  considered in the first barycentric subdivision of $K(\Gamma)$.}  $n$-manifold iff, for each color $c\in\Delta_n$, all $\hat c$-residues of $\G$ represent closed connected $(n-1)$-manifolds.
\end{itemize}

\begin{remark} \label{correspondence-sing-boundary}{\em  If $N$ is a singular $n$-manifold, then a compact $n$-manifold $\check N$ is easily obtained by deleting small open neighbourhoods of its singular vertices.
Conversely, given a compact $n$-manifold $M$, a singular $n$-manifold $\widehat M$ can be constructed by capping off each component of $\partial M$ by a cone over it.

Throughout the paper, we will restrict to the class of compact $4$-manifolds without spherical boundary components, where the above correspondence is bijective; hence, singular $n$-manifolds and compact $n$-manifolds can be associated to each other in a well-defined way. 
Obviously, in this wider context, closed $n$-manifolds are characterized by $M= \widehat M.$  }
\end{remark}

In virtue of the bijection described in Remark \ref{correspondence-sing-boundary}, an $(n+1)$-colored graph $\G$ is said to {\it represent}
a compact $n$-manifold $M$ 
(or, equivalently, to be a {\it gem} of $M$, where gem means {\it Graph Encoding Manifold})  if and only if  it represents the associated singular manifold $\widehat M$.  

A {\it $\hat c$-residue} of a gem $\G$ will be called {\it singular} if it corresponds to a singular vertex of $\vert K(\G) \vert$, i.e. if it does not represent the sphere, while a color $c$ is said to be {\it singular} if at least one $\hat c$-residue is singular.  

\medskip

{\it Gem theory} (or {\it crystallization theory}) is based on the following existence theorem, which extends to the boundary case the foundational one due to Pezzana for closed manifolds of arbitrary dimension (see \cite{Ferri-Gagliardi-Grasselli}).

\begin{theorem}{\em (\cite{Casali-Cristofori-Grasselli})}\ \label{Theorem_gem}  
Any compact $n$-manifold $M$ 
admits a gem $\G$, that is bipartite if and only if $M$ is orientable.
\par \noindent In particular, if  $M$ has empty or connected boundary: 
\begin{itemize}
\item $\G$ may be assumed to have color $n$ as its unique possible singular color, and exactly one $\hat n$-residue (in this case we will say that $\G$ belongs to the class $G^{(n)}_s$);
\item $\G$ may be assumed to have exactly one $\hat c$-residue, $\forall c \in \Delta_n$ (in this case $\G$ is  called  
a \emph{crystallization} of $M$).
\end{itemize}
\end{theorem}

\medskip

Within gem theory, a finite set of combinatorial moves have been defined, which translate the homeomorphism problem of the represented manifolds. 

\begin{definition}\label{def_dipole}{\em An {\it $r$-dipole ($1\le r\le n$) of colors $c_1,\ldots,c_r$} of an $(n+1)$-colored graph $\G$ is a subgraph consisting of two vertices, belonging to different connected components of $\Gamma_{\hat c_1\ldots\hat c_r}$, that are joined by $r$ edges, colored by $c_1,\ldots,c_r.$

The {\it elimination} of an $r$-dipole in $\Gamma$ can be carried out by deleting the subgraph and welding the remaining hanging edges according to their colors; in this way another $(n+1)$-colored graph $\Gamma^\prime$ is obtained. The inverse operation is called the {\it addition} of the dipole to $\Gamma^\prime.$}
\end{definition}

\begin{proposition}
Let $\G$ be a gem of a compact $n$-manifold $M$ with empty or connected boundary.  If $\Gamma^\prime$ is obtained from $\Gamma$ by eliminating an $r$ dipole ($1\le r\le n$), then $\Gamma^\prime$ is a gem of $M$, too. 
\end{proposition}

Moreover, other combinatorial moves have been defined, whose effects on the represented manifolds turn out to be significant also in the context of the present paper (see Proposition \ref{chiudere}, which gives an original contribution to the study of these moves, applied to colored graphs representing compact manifolds with non-empty boundary). 

\begin{definition}\label{rho-pair}{\em A $\rho_h${\it -pair} ($1\leq h\leq n$) of color $i\in\Delta_n$ in an $(n+1)$-colored graph $\Gamma$ is a pair of $i$-colored edges $(e,f)$ sharing the  same $\{i,c\}$-colored cycle for each $c\in\{c_1,\ldots,c_h\}\subseteq\Delta_n.$ 
\\
The {\it switching} of $(e,f)$ consists in canceling $e$ and $f$ and establishing new $i$-colored edges between their endpoints; the reversed operation is obviously the switching of a $\rho_{n-h}$-pair.  
Although, in general, the switching can 
be performed in two different ways, it is uniquely determined if $\G \in G^{(n)}_s,\ h\in \{n-1, n\}$ and the bipartition of each non-singular $\hat c$-residue is preserved.  }
\end{definition} 

The topological effects of the switching of $\rho_{n-1}$- and $\rho_n$-pairs have been completely determined in the case of closed $n$-manifolds: see \cite{Bandieri-Gagliardi}, where it is proved that a $\rho_{n-1}$-pair (resp. $\rho_n$-pair) switching does not affect the represented $n$-manifold  (resp. either induces the splitting into two connected summands, or the ``loss'' of a $\mathbb S^{1} \otimes \mathbb S^{n-1}$  summand in the represented $n$-manifold). 

\bigskip 

Finally, we point out that an interesting aspect of gem theory is related to the possibility of combinatorially defining PL invariants of manifolds in arbitrary dimension.  The most important one, which extends to higher dimension the classical genus of a surface and the Heegaard genus of a $3$-manifold, is founded on the existence of a particular type of embedding of colored graphs into surfaces. 

\begin{proposition}{\em (\cite{Ferri-Gagliardi-Grasselli})}\label{reg_emb}
Let $\G$ be a connected bipartite (resp. non-bipartite) 
$(n+1)$-colored graph of order $2p$. Then for each cyclic permutation $\varepsilon = (\varepsilon_0,\ldots,\varepsilon_n)$ of $\Delta_n$, up to inverse, there exists a cellular embedding, called \emph{regular}, of $\G$  
into an orientable (resp. non-orientable)
closed surface $F_{\varepsilon}(\G)$ whose regions are bounded by the images of the $\{\varepsilon_j,\varepsilon_{j+1}\}$-colored cycles, for each $j \in \mathbb Z_{n+1}$.
Moreover, the genus (resp. half the genus)  of $F_{\varepsilon}(\G)$, denoted by $\rho_{\varepsilon} (\G)$,  satisfies

\begin{equation*}
2 - 2\rho_\varepsilon(\G)= \sum_{j\in \mathbb{Z}_{n+1}} g_{\varepsilon_j, \varepsilon_{j+1}} + (1-n)p.
\end{equation*}

\end{proposition}

Throughout the paper, for sake of notational simplicity, if $i \in \Delta_n$ is such that $g_{\hat i}=1$, we will write $\rho_{\varepsilon}(\G_{\hat\imath})$  (resp. $F_{\varepsilon}(\G_{\hat\imath})$) instead of  $\rho_{\varepsilon^\prime}(\G_{\hat\imath})$ (resp. $F_{\varepsilon^\prime}(\G_{\hat\imath})$), $\e^\prime$ being the permutation induced by $\e$ on $\hat\imath$. 

\medskip 

\begin{definition} {\em The \emph{regular genus} of  an $(n+1)$-colored graph $\G$ is defined as
$$\rho(\G) = min\{\rho_\varepsilon(\G)\ \vert \ \varepsilon\ \text{cyclic permutation of \ } \Delta_n\},$$
while the  {\it regular genus} of a compact $n$-manifold $M$ is defined as}
$$\mathcal G (M) = min\{\rho(\G)\ \vert \ \G\ \text{gem of \ } M\}.$$
\end{definition}

\bigskip 

Since the number of vertices of a gem, i.e. the number of $n$-simplices in the dual triangulation, can be interpreted as a measure of the ``complexity" of the represented manifold, another PL invariant naturally arises:

\begin{definition} \label{def. gem-complexity}
{\em For each compact $n$-manifold $M$, its \emph{gem-complexity} is the non-negative integer $k(M)= p - 1$, 
where $2p$ is the minimum order of an $(n+1)$-colored graph representing $M$.}
\end{definition}


\begin{remark} \label{rem. invariants}{\em  Many significant classification results have been obtained within crystallization theory,  with respect to both regular genus and gem-complexity:  see, for example, \cite{Casali-Cristofori-Gagliardi Complutense 2015} for a review regarding  dimension $4$. Regular genus and gem-complexity are also proved to be related to the {\it G-degree}, a PL invariant arising within the theory of {\it Colored tensor models} in theoretical physics: this relation is investigated - among other - in \cite{Casali-Cristofori-Dartois-Grasselli} and \cite{Casali-Cristofori-Grasselli}.   
The topic of the present paper highlights further connections of the above invariants with those based on other representation tools for PL $4$-manifolds.   }
\end{remark}


\section{Relationships and results}


\subsection{From Kirby diagrams to gems}  \label{ss.from_Kirby_to_gems}
\medskip

In this section we will describe an algorithm to obtain a gem $\G(L,d)$ of a compact (orientable) 4-manifold $M$ starting from a Kirby diagram $(L,d)$ of $M$, where $L$ is a link in $\mathbb S^3$ and $d$ a vector that contains the framings. More precisely, we suppose that $L$ has $l$ components, $L_1,\ldots,L_l$, whose first $m$, if any, are dotted and, for each $i\in\{1,\ldots,l-m\}$, $d_i$ is the framing of $(m+i)$-th component.

Actually, the construction is performed on a planar diagram of $L$, which we suppose fixed and, for sake of simplicity, will also be called $L$; hence we will refer to all combinatorial characteristics of the diagram (such as {\it arcs, crossings} and {\it regions}) as of $L$ itself.  
Moreover,  we assume that the diagram is connected\footnote{We point out that any Kirby diagram can be made so by suitable Reidemeister moves. Alternatively, the algorithm can be applied to each connected component and then graph-connected sums can be performed on the resulting gems in order to obtain a 5-colored graph representing $M$ (see \cite{Casali-Cristofori trisection bis} for details).} and {\it self-framed}, i.e. the writhe of each non-dotted component equals its framing. Note that this last requirement can be always satisfied up to the addition of positive or negative curls. 

\begin{figure} [h!]
    \centering
    \includegraphics[width=0.4\linewidth]{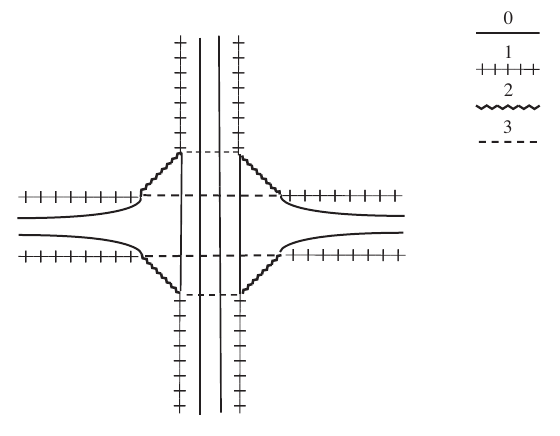}
\caption{subgraph corresponding to a crossing}
    \label{fig:crossing}
\end{figure}

\begin{figure} [h!]
    \centering
    \includegraphics[width=0.95\linewidth]{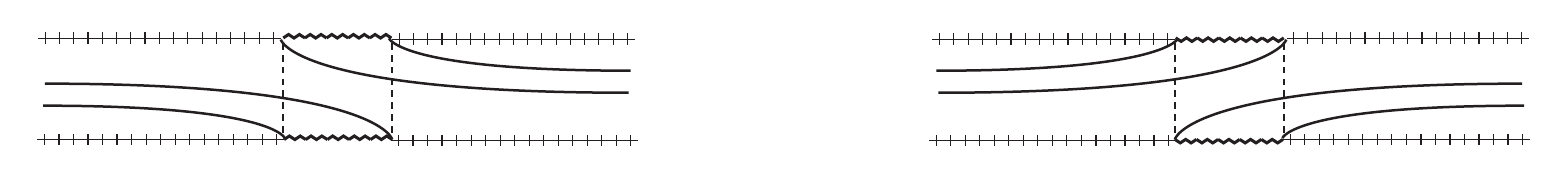}
\caption{subgraphs corresponding to a positive (left) or negative (right) curl}
    \label{fig:positive-negative-graph-curl}
\end{figure}

\begin{figure} [h!]
    \centering
    \includegraphics[width=0.2\linewidth]{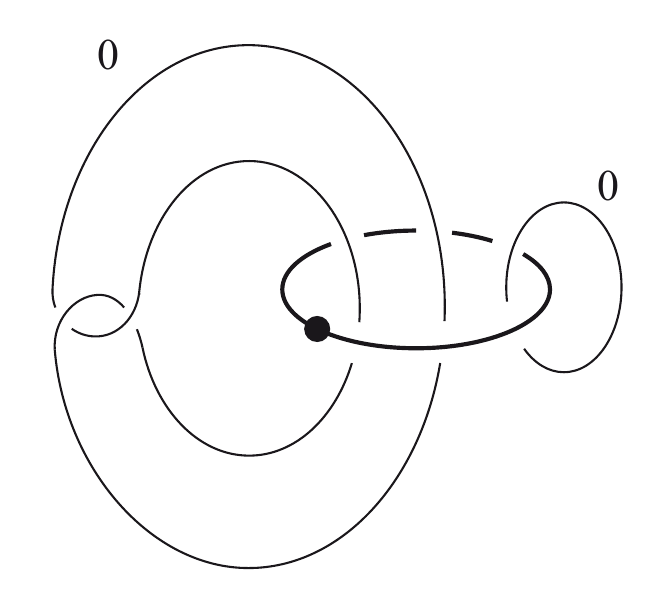}\quad
    \includegraphics[width=0.7\linewidth]{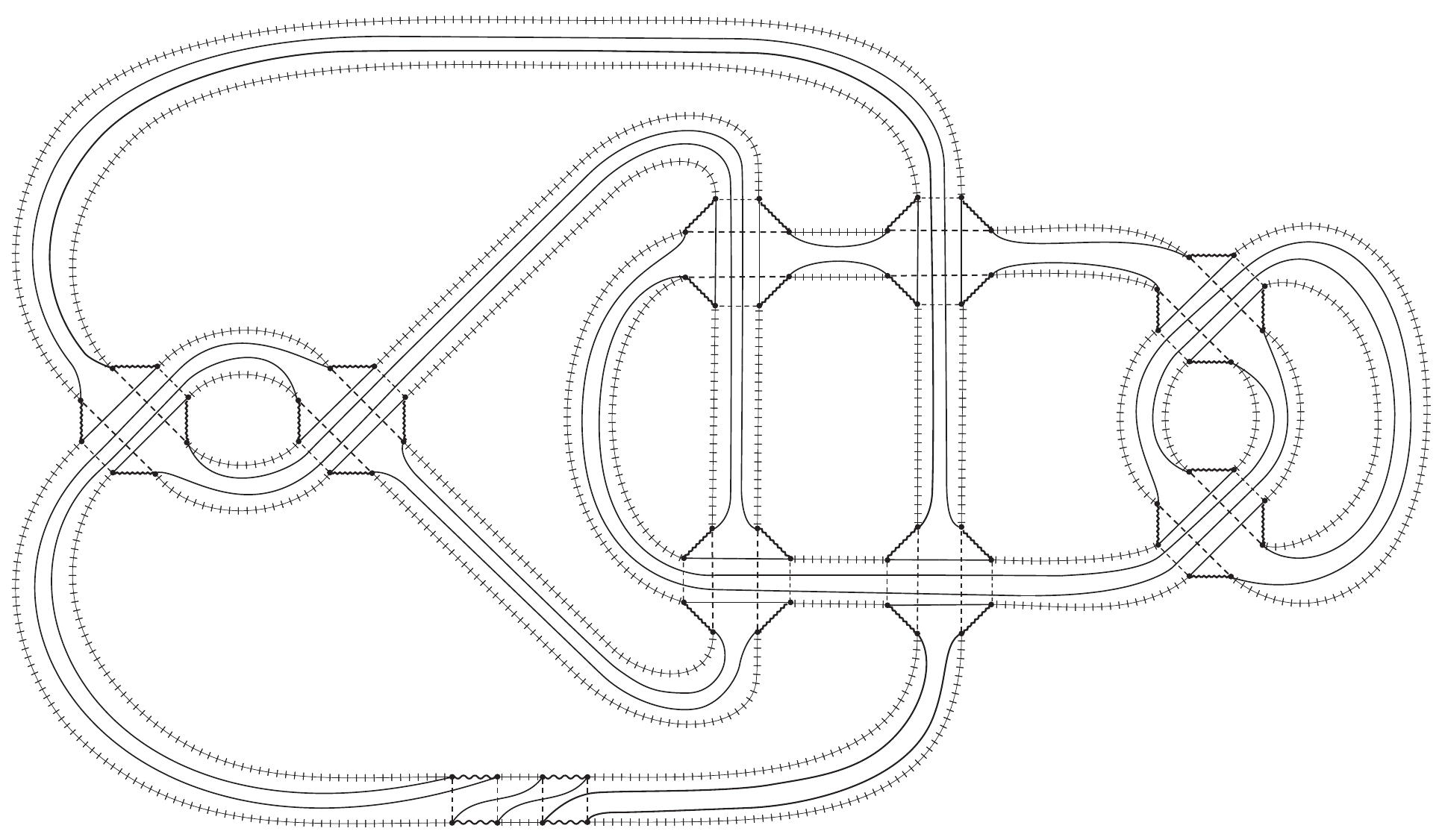}
\caption{A Kirby diagram of $\mathbb S^2\times\mathbb D^2$ and a gem of its boundary}
    \label{fig:fishtail_3dim}
\end{figure}

The starting step of the procedure is the construction of a 4-colored graph $\Lambda(L,d)$ representing the 3-manifold $\partial M$, which, as recalled in Section \ref{ss:Kirby diagrams}, is obtained by Dehn surgery on $\mathbb S^3$ along the framed link associated to $L.$ As it is easily seen from the example in Figure \ref{fig:fishtail_3dim}, $\Lambda(L,d)$ can be directly ``drawn over" $L$, by taking the subgraph of order eight depicted in Figure \ref{fig:crossing} for each crossing, while to each curl there corresponds a subgraph of one of the two types in Figure \ref{fig:positive-negative-graph-curl} according to the curl being positive or negative.

A central role is played in the algorithm by a combinatorial structure, called a {\it quadricolor}, consisting of four vertices $\{P_0,\ P_1,\ P_2,\ P_3\}$ such that, for each $i\in\mathbb Z_4$, $P_{i}$ and $P_{i+1}$ are $i$-adjacent and $P_{i}$ does not belong to the $\{i+1,i+2\}$-colored cycle shared by the other three vertices. We point out that a quadricolor always arises between the subgraphs corresponding to a curl and an undercrossing or between two curls; therefore, up to the addition of a pair of opposite curls, the existence of a quadricolor can always be assumed on each framed component of $L.$

\begin{figure} [h!]
    \centering
        \includegraphics[width=0.9\linewidth]{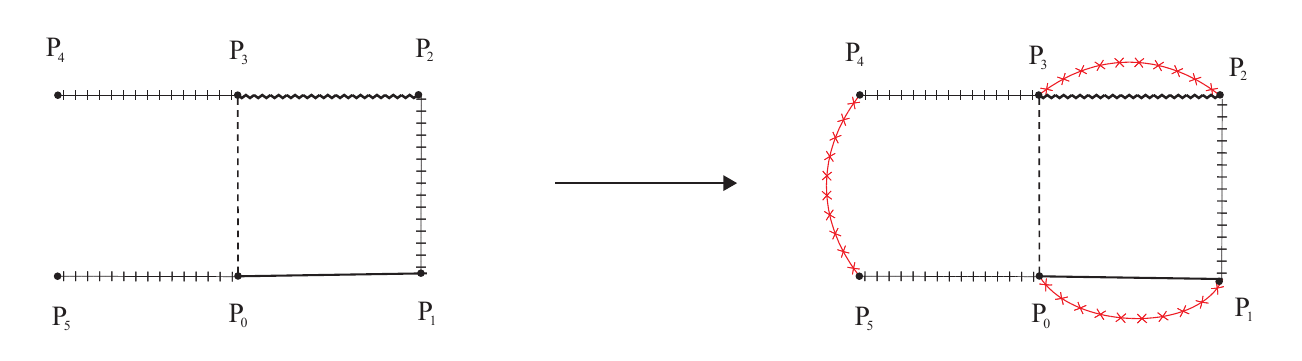}
\caption{The attachment of 4-colored edges among the vertices of a quadricolor}
    \label{fig:quadricolor}
\end{figure}

\begin{figure} [h!]
    \centering
        \includegraphics[width=0.7\linewidth]{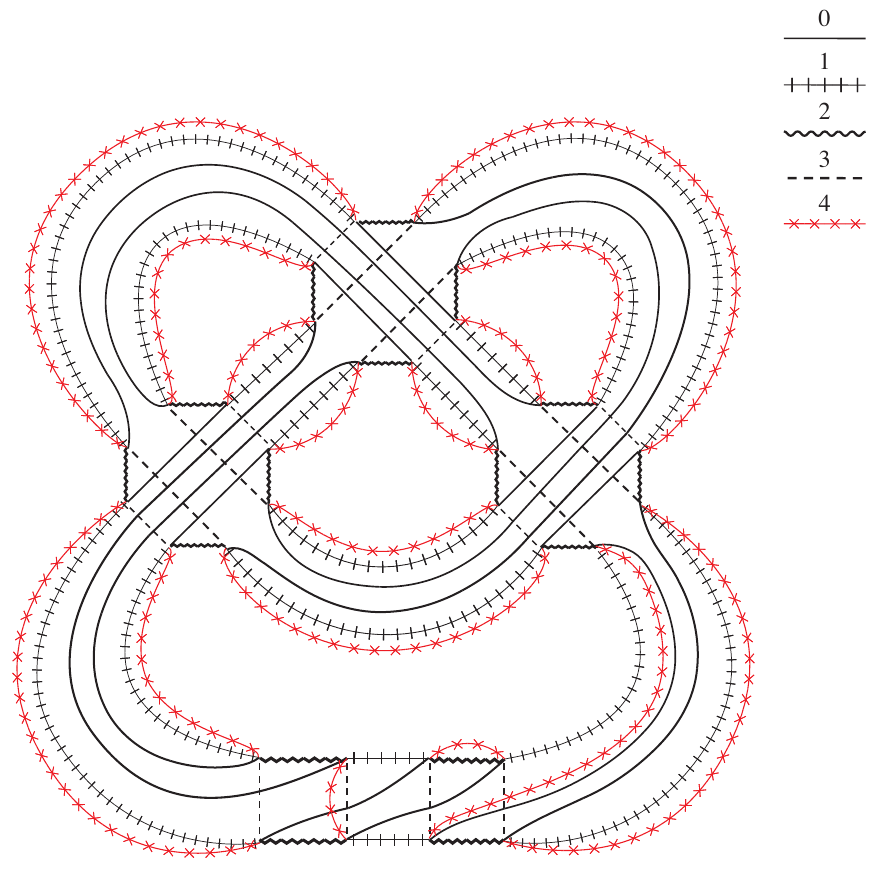}  
\caption{The 5-colored graph $\G(L,d)$ obtained from the framed link $({K_T},+1)$, $K_T$ being the trefoil knot}  
    \label{fig:M4_trefoil1}
\end{figure}

In the special case of $m=0$, i.e. when dotted components are missing, the procedure to obtain a gem $\G(L,d)$ of $M$ is particularly simple: it is sufficient to double by 4-colored edges all 1-colored edges of $\Lambda(L,d)$ except among the vertices of a quadricolor for each (framed) component, where the addition of 4-colored edges is done as shown in Figure \ref{fig:quadricolor} (see also Figure \ref{fig:M4_trefoil1} for an example where $(L,d)$ is a trefoil with framing +1).

The case when dotted components appear requires some further hypotheses and some choices of points and parts of arcs of the diagram, each choice producing at the end a different graph but all representing the compact manifold encoded by the Kirby diagram.

The dotted components are assumed to be in {\it good position}, i.e. they are unknotted, unlinked and on each of them overcrossings and undercrossings never alternate\footnote{This assumption does not involve any loss of generality since, by suitably arranging handles and possibly using Reidemeister moves, any compact orientable 4-manifold turns out to admit a Kirby diagram whose dotted components are in good position.}.

\begin{figure} [h!]
    \centering
    \includegraphics[width=0.35\linewidth]{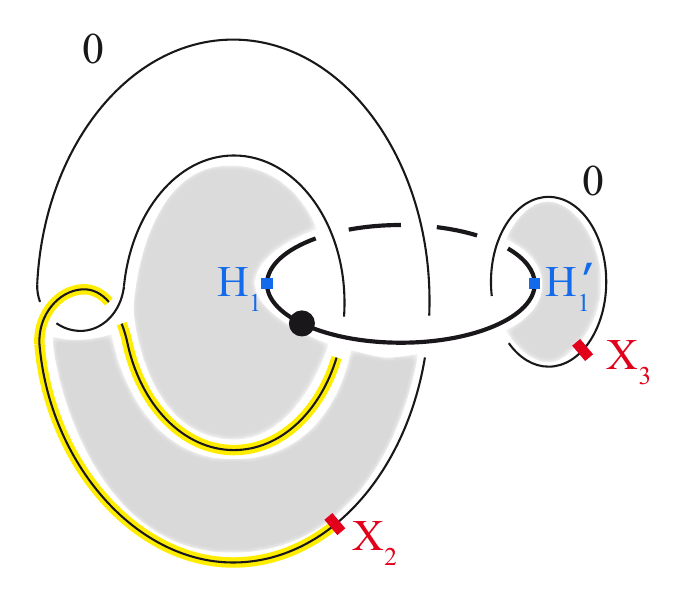}
\caption{}
    \label{fig:fishtail_choices}
\end{figure}

\begin{figure} [h!]
    \centering
    \includegraphics[width=0.47\linewidth]{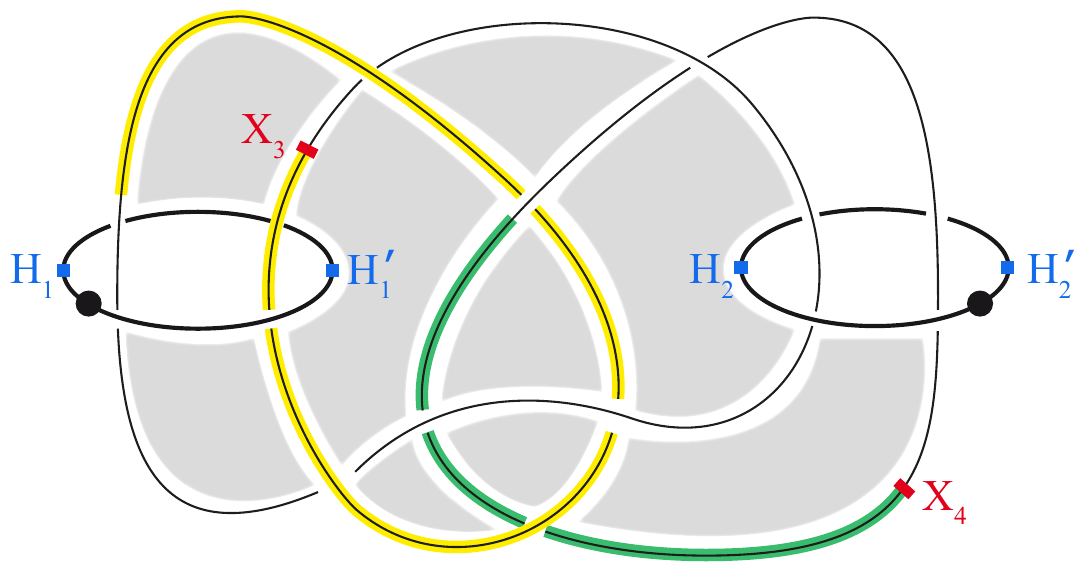}\quad\includegraphics[width=0.47\linewidth]{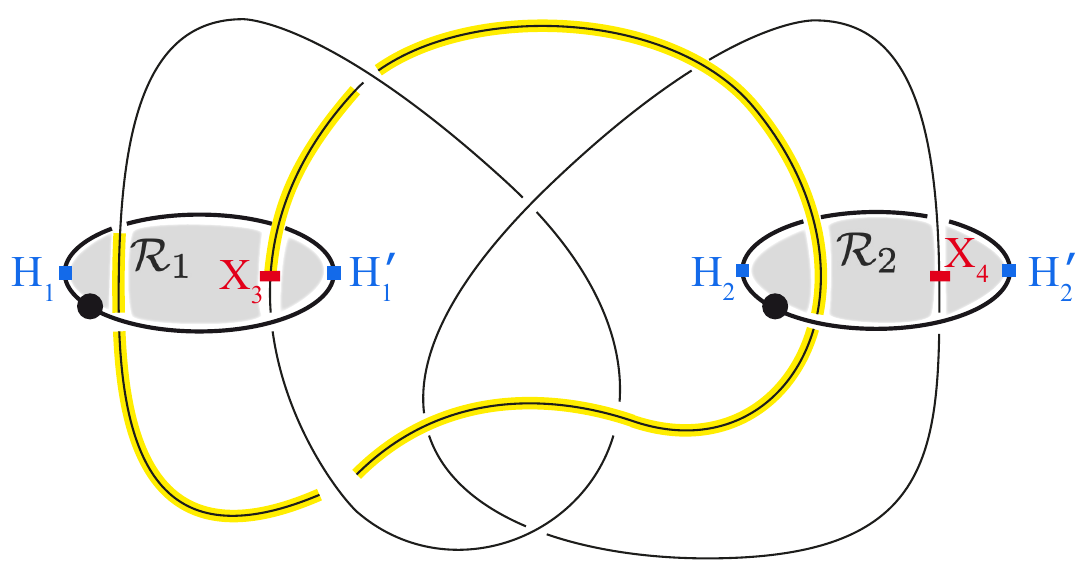}
\caption{Different choices on the same diagram: $\mathcal R_1=\mathcal R_2$ on the left; $\mathcal R_1\neq\mathcal R_2$ and an empty sequence of highlighted segments for $L_4$ on the right}
    \label{fig:K0&trefoil}
\end{figure}

The choices we are listing below will determine how to add the 4-colored edges to the gem $\Lambda(L,d)$  of $\partial M$ which has been previously constructed:
\begin{itemize}
\item on each dotted component $L_i$ ($i\in\{1,\ldots,m\}$) two points $H_i$ and $H'_i$ must be ``marked", such that they divide it into two parts, one containing only overcrossings and the other containing only undercrossings;
\item for each $i\in\{m+1,\ldots,l\}$ the framed component $L_i$ is ``cut" in a point $X_i$ between a curl and an undercrossing or between two curls; moreover, starting from $X_i$ in the direction opposite to the undercrossing (or in any direction if $X_i$ lies between two curls) let’s ``highlight" a sequence of consecutive segments of arcs of $L_i\setminus X_i$, so that for each dotted component $L_j$ the points $H_j$ and $H'_j$ belong to the boundary of the same region $\mathcal R_j$ of the plane determined by $L$ with all cuts and highlighted segments removed. 
\end{itemize}
Figure \ref{fig:fishtail_choices} shows an example of choices where the determined region, $\mathcal R_1$, is the union of the gray ones with the unbounded region. Figure \ref{fig:K0&trefoil} shows examples of different choices on the same Kirby diagram.

The above choices allow to determine the attachments of the 4-colored edges to $\Lambda(L,d)$ in order to obtain $\G(L,d)$ as follows:
\begin{itemize}
\item for each $i\in\{m+1,\ldots,l\}$, a quadricolor can be detected in $\Lambda(L,d)$ in correspondence with the cut $X_i$ and so the triple of 4-colored edges among its vertices is added according to Figure \ref{fig:quadricolor};
\item for each quadricolor, starting from the two vertices 1-adjacent to $P_4$ and $P_5$ (see Figure \ref{fig:quadricolor}), 4-colored edges are added, following the sequence of highlighted segments of arcs and respecting bipartition classes, among the vertices of each pair of 1-colored edges which are ``parallel" to a highlighted segment; 
\item the 0-colored edges corresponding to any undercrossing which is ``met" by the sequence are doubled by 4-colored edges;
\item for each $i\in\{1,\ldots,m\}$ the points $H_i$ and $H_i^\prime$ allow to identify, on the subgraph corresponding to the dotted component $L_i$, two pairs of vertices $(v_i,w_i)$ and $(v_i',w_i')$, which are the endpoints of the two 1-colored edges separating overcrossings from undercrossings and lie on the side of the region $\mathcal R_i$; then a 4-colored edge is added between $v_i$ and $v_i'$ and between $w_i$ and $w_i'$; 
\item the missing 4-colored edges are added so as to obtain a cycle from each $\{1,4\}$-colored path.
\end{itemize}

\begin{figure} [h!]
    \centering
    \includegraphics[width=0.9\linewidth]{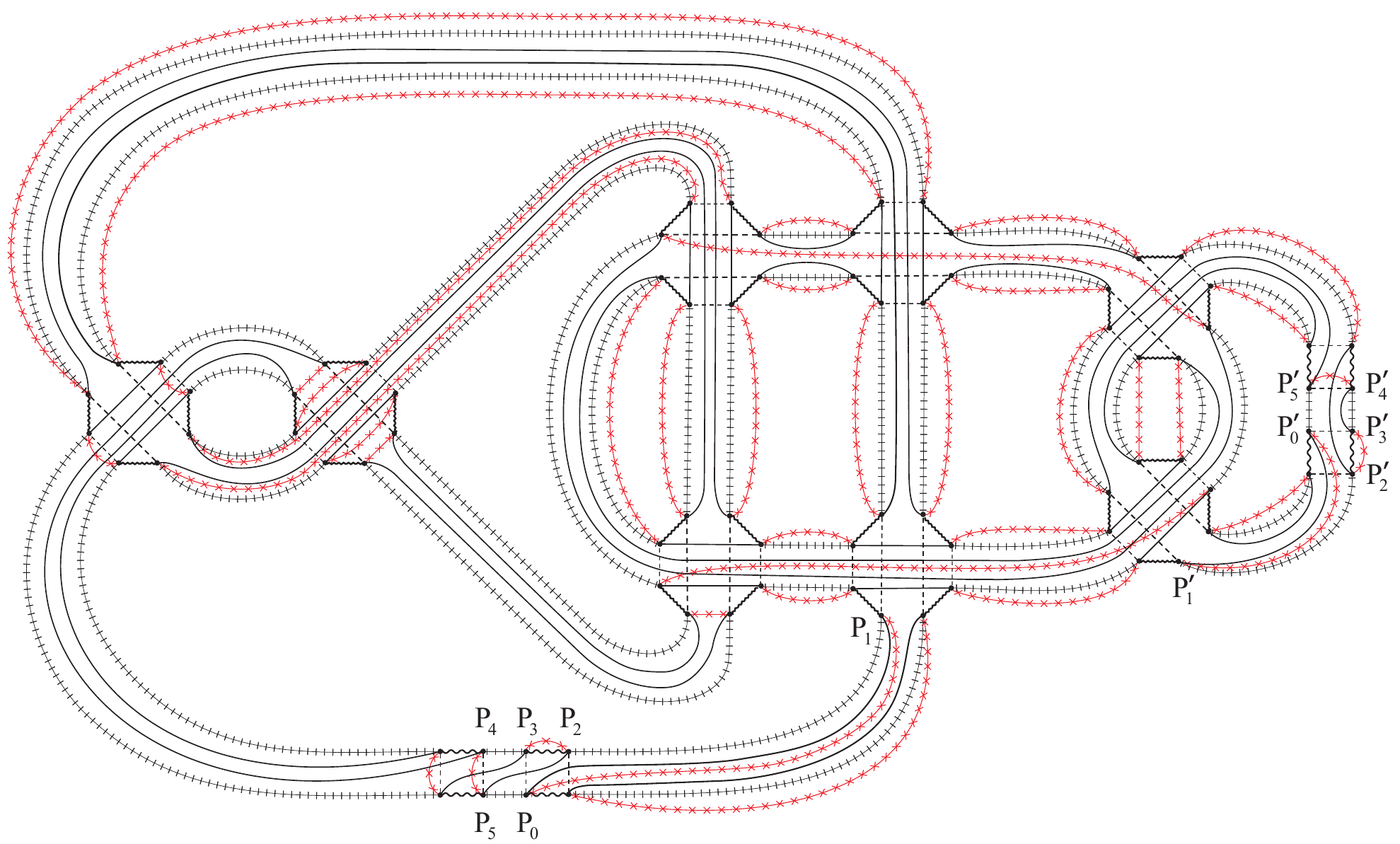}
\caption{}
    \label{fig:fishtail_4dim}
\end{figure}

Figure \ref{fig:fishtail_4dim} shows the application of the above algorithm to the example of Figure \ref{fig:fishtail_choices}. 

\begin{theorem} \ {\rm (\cite{Casali-Cristofori Kirby-diagrams})} \ \label{t.Kirby-diagram}
Given a compact 4-manifold $M$, let $(L,d)$ be a connected Kirby diagram representing it, where all dotted components, if any, are in good position. Then, the 5-colored graph $\G(L,d)$ is a gem of $M.$  
\end{theorem}

\begin{figure} [h!]
    \centering
    \includegraphics[width=0.9\linewidth]{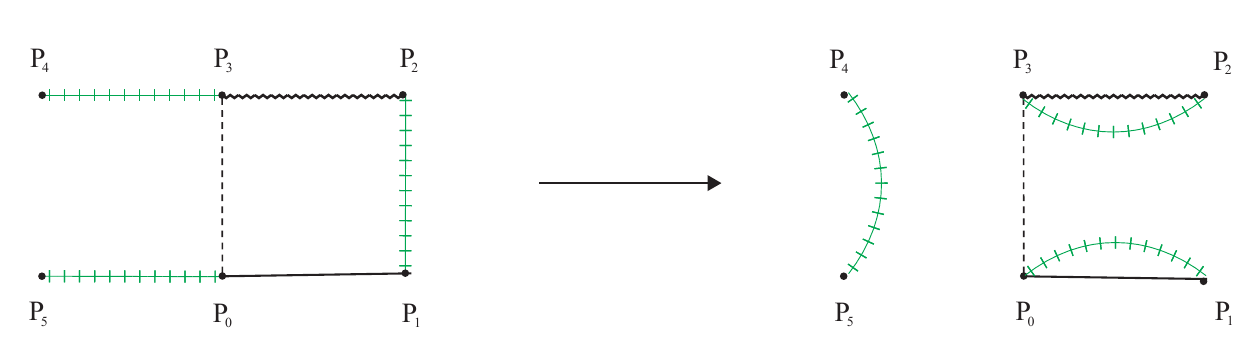}
\caption{the smoothing of a quadricolor}
    \label{fig:smoothing}
\end{figure}

The proof of the above result consists in showing that $\G(L,d)$ can also be obtained from $\Lambda(L,d)$ by a sequence of combinatorial moves with a precise topological meaning.
One of the most important is the so-called {\it smoothing of a quadricolor} which is the exchange of three 1-colored edges among the vertices of the quadricolor as depicted in Figure \ref{fig:smoothing}.

This move, done on a quadricolor of $\Lambda(L,d)$ corresponding to a given framed component $L_i$, produces a new 
4-colored graph representing the 3-manifold obtained by Dehn surgery on the link obtained from $L$ by adding the complementary knot of $L_i$ (i.e. a $0$-framed trivial knot linking the component $L_i$ geometrically once), that is, essentially, by cancelling that framed component\footnote{Quadricolors in 4-colored graphs were originally introduced by Lins (\cite{Lins-book}). The smoothing of a quadricolor in a gem of a closed 3-manifold can also be interpreted as the substitution, in the dual triangulation, of a solid torus with another solid torus having the same boundary. Hence, it is equivalent to performing a Dehn surgery. In the specific case of $\Lambda(L,d)$, this surgery can be identified as the one along the complementary knot of the component naturally associated with the quadricolor (see \cite[Proposition 9(i)]{Casali-Cristofori Kirby-diagrams} for details).}.  

\bigskip

{\it Sketch of the proof (of Theorem \ref{t.Kirby-diagram})}:

First, let us note that, as visualized in Figure \ref{fig:schizzo-dim}, the handle-decomposition encoded by $(L,d)$ can be realized by considering a collar of $\#_m (\mathbb S^1\times\mathbb S^2)$, where $m\geq 0$ is the number of dotted components of $(L,d)$, and attaching the $0$- and $1$-handles on one of its boundary component, the $2$-handles on the other one.  
This can be done via gems in the following way.

\begin{figure} [h!]
    \centering
    \includegraphics[width=1\linewidth]{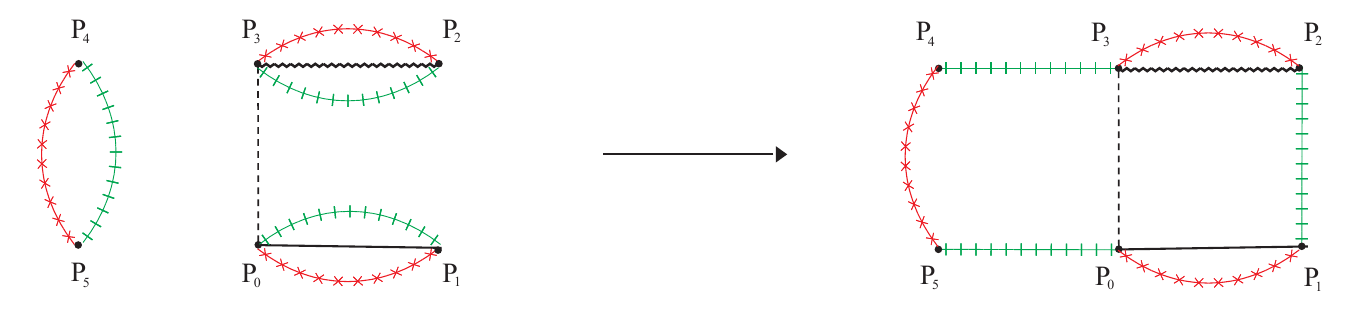}
\caption{the attachment of a 2-handle}
    \label{fig:2-handle}
\end{figure}

\begin{figure}  [h!]
    \centering
    \includegraphics[width=0.9\linewidth]{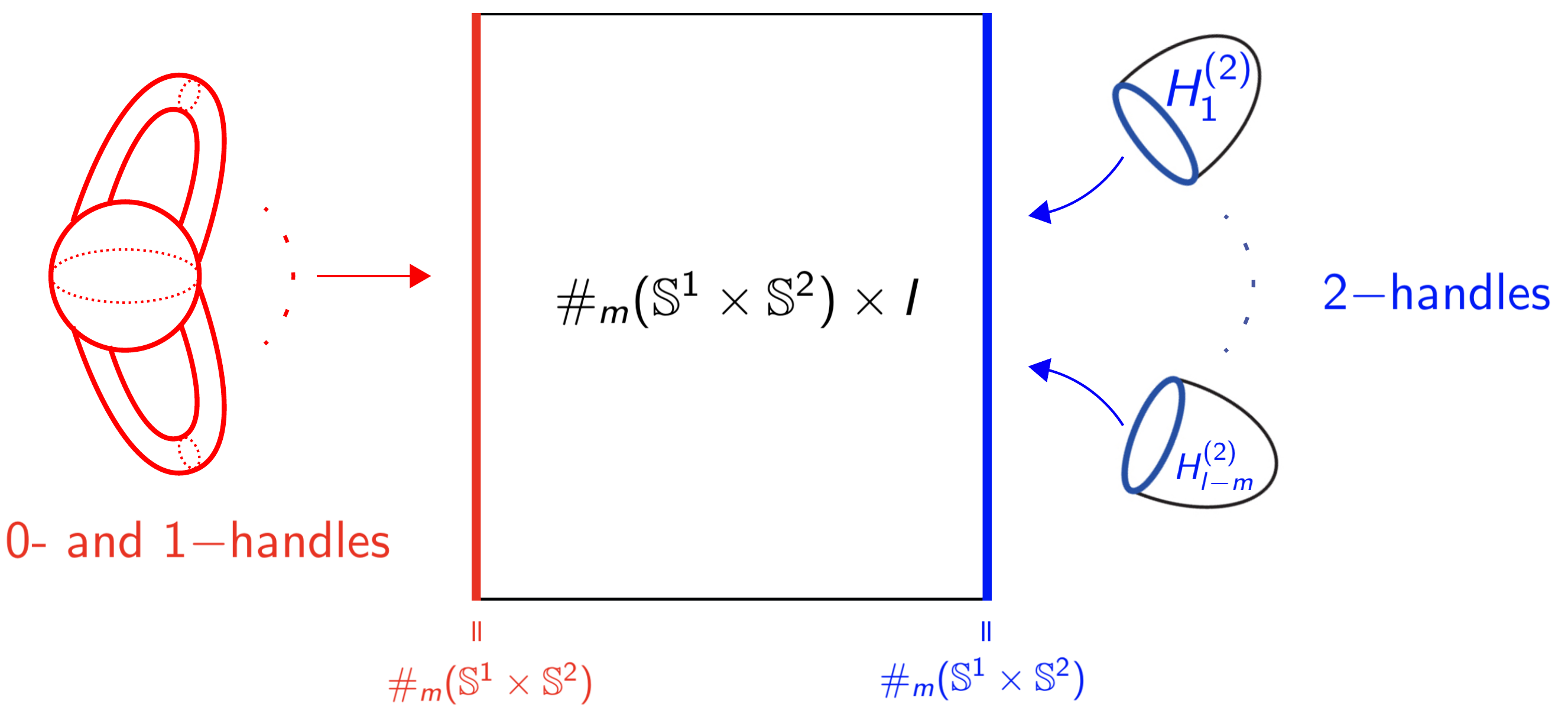}
    \caption{Realizing the handle-decomposition encoded by a Kirby diagram via gems (see proof of Theorem \ref{t.Kirby-diagram})}
    \label{fig:schizzo-dim}
\end{figure}

\begin{itemize}
    \item Starting from $\Lambda(L,d)$, let us perform the smoothing of all quadricolors of the framed components of $(L,d)$; by what we have just observed, the result is a 4-colored graph $\tilde\Lambda$ representing $\#_m (\mathbb S^1\times\mathbb S^2).$ 
By doubling the 1-colored edges of $\tilde\Lambda$ by 4-colored ones, a 5-colored graph $\tilde\G$ is obtained which can be easily seen to represent $\#_m (\mathbb S^1\times\mathbb S^2)\times I$, whose two boundary components are represented by the 4-colored graphs $\tilde\G_{\hat 1}$ and $\tilde\G_{\hat 4}$ respectively (see the central part of Figure \ref{fig:schizzo-dim}).

\item  
A sequence of $\rho_2$-pairs of color 4 (see Section \ref{preliminaries.gems}) can be detected in $\tilde\G$, whose switchings are proved to preserve the represented manifold.  
In fact, the resulting 5-colored graph can also be obtained by {\it suturing} a sequence of {\it wounds} in $\tilde \Lambda$, thought of as a 5-coloured graph with boundary, and subsequently by {\it capping it off with respect to color 1} 
(see  \cite{Gagliardi 1987} and \cite{generalized-genus} for a detailed description of the cited moves). The fact that these operations do not change the PL type of the represented 4-manifold (they simply produce a re-triangulation of it) is mainly due to the existence in $\tilde\Lambda$ of a standard sequence of dipole eliminations corresponding to the highlighted segments of arcs (see \cite[Remark 5]{Casali-Cristofori Kirby-diagrams}).

\item As a consequence of the $\rho_2$-pairs switchings, in each subgraph corresponding to a dotted component $L_i$ ($i\in\{1,\ldots,m\}$) a $\rho_3$-pair of color 4 appears, whose endpoints are exactly $(v_i,w_i)$ and $(v_i',w_i')$.
As pointed out in \cite[Proposition 11]{Casali-Cristofori Kirby-diagrams}, the switching of these $\rho_3$-pairs, realized by connecting $v_i$ with $v_i'$ and $w_i$ with $w_i'$ by 4-colored edges for each $i\in\{1,\ldots,m\}$, has the topological effect of an identification between the boundary of a genus $m$ 4-dimensional handlebody (that is the union of 0- and 1-handles in the handle-decomposition encoded by $(L,d)$) and the component of the boundary of $\#_m (\mathbb S^1\times\mathbb S^2)\times I$ represented by $\tilde\G_{\hat 1}$ (the ``red" boundary in Figure \ref{fig:schizzo-dim}) after the $\rho_2$-pairs switchings.  

\item  Note that $\tilde\G_{\hat 4}$, the ``blue" boundary in Figure \ref{fig:schizzo-dim}, has not been affected by the performed switchings of $\rho_2$- and $\rho_3$-pairs and, therefore, can be thought of as representing the boundary of the resulting 4-manifold (which is obviously homeomorphic to a genus $m$ 4-dimensional handlebody). 

\noindent In order to obtain $\G(L,d)$, it is now sufficient to re-place the triples of 1-colored edges of the quadricolors in their original position (see Figure \ref{fig:2-handle}). 
As proved in \cite[Proposition 9(ii)]{Casali-Cristofori Kirby-diagrams}, for each framed component $L_i$ ($i\in\{m+1,\ldots,l\}$), this move on the associated quadricolor realizes the attachment of a 2-handle along $L_i$ on the boundary component represented by $\tilde\G_{\hat 4}.$
\end{itemize}
\ \qed

As a direct consequence of the construction, the order of $\G(L,d)$ and its regular genus can be expressed in terms of the combinatorial properties of the Kirby diagram. Therefore, upper bounds for the value of the gem-complexity and the regular genus of the represented compact 4-manifold can be obtained:

\begin{theorem} \ {\rm (\cite{Casali-Cristofori Kirby-diagrams})} \  \label{regular-genus&gem-complexity(dotted)} 
Let $M$ be a compact 4-manifold and let $(L,d)$ be a connected Kirby diagram representing $M$, with $l$ components whose dotted ones, if any, are in good position. Then:
$$\mathcal G(M)\leq  \  s + \bar s  + (l-m) +1$$ 
\noindent where $m \ge 0$ is the number of dotted components, $s$ the number of crossings of $(L,d)$ and $\bar s$ the number of undercrossings of its framed components.
Moreover, if $L$ is different from the trivial knot: 
$$k(M)\leq \ 2s + 2 \bar s  + 2 m -1 + 2\sum_{i=1}^{l-m} t_i  $$
\noindent where 
$t_i=2$ if the writhe and the framing of the $i$-th framed component coincide, otherwise $t_i$ is the absolute value of their difference.
\end{theorem} 

We 
point out that, as proved in \cite{Casali JKTR2000}, if $(L,d)$ has no dotted components, the 4-colored graph $\Lambda(L,d)$ contains some particular combinatorial configurations, called {\it generalized dipoles} that can be manipulated, thus obtaining a new 4-colored graph $\Omega(L,d)$, without changing the represented 3-manifold or affecting the quadricolors; by doubling by 4-colored edges the 1-colored ones and performing the moves on quadricolors of $\Omega(L,d)$ (i.e.  by applying the construction to $\Omega(L,d)$ instead of $\Lambda(L,d)$), a 5-colored graph $\tilde\Omega(L,d)$ is obtained, still representing $M$, whose combinatorial structure allows to establish the following improvement of the estimation of the regular genus:

\begin{theorem} \ {\rm (\cite{Casali-Cristofori Kirby-diagrams})} \  \label{regular-genus&gem-complexity(framed)} 
Let $M$ be a compact 4-manifold and let $(L,d)$ be a connected Kirby diagram  with no dotted components representing $M$. Then:
$$\mathcal G(M)\leq m_{\alpha} + l$$
\noindent where $l$ is the number of components of $L$ and $m_\alpha$ is the number of $\alpha$-colored regions in a chess-board coloration of $L$, $\alpha$ being the color of the unbounded region.
\end{theorem}

The following result - which can be easily proved by means of the above constructions (for details, see the proofs of \cite[Theorem 1]{Casali JKTR2000}, \cite[Theorem 7]{Casali-Cristofori Kirby-diagrams} or \cite[Theorem 3]{Casali-Cristofori gem-induced}) -  will turn out to be useful in Section \ref{ss.from_Kirby_to_trisections}. 

\begin{proposition}\label{order-regular genus Gamma}\ Let $(L,d)$ be a connected Kirby diagram, with $s$ crossings, whose dotted components, if any, are in good position. If the cyclic permutation $\varepsilon = (1,0,2,3,4)$ is chosen, then:
$$\rho_{\varepsilon}(\G_{\hat 4}(L,d))=\rho_{\varepsilon}(\Lambda(L,d)) = s + 1.$$
\smallskip
Moreover, if dotted components are missing and $m_\alpha$ is the number of $\alpha$-colored regions in a chess-board coloration of $L$ (where $\alpha$ is the color of the unbounded region), then $$\rho_{\varepsilon}(\tilde\Omega_{\hat 4}(L,d))=\rho_{\varepsilon}(\Omega(L,d)) = \  m_\alpha.$$
\end{proposition}

\bigskip
\bigskip

As already pointed out in Section \ref{ss:Kirby diagrams}, if $\partial M\cong\#_r(\mathbb S^1\times\mathbb S^2)$, a Kirby diagram $(L,d)$ of $M$ also represents the closed 4-manifold $\bar M$ uniquely obtained by adding $r$ 3-handles and a 4-handle to the handle decomposition encoded by $(L,d)$. 

Note that, if $r=0$, then $\G(L,d)$ is also a gem of $\bar M$ (identified with $M$ by capping off the spherical boundary). On the other hand, if $r>0$, a gem of $\bar M$ can be directly obtained from $\G(L,d)$ in some particular cases, namely, for example, if $\G(L,d)$ contains $r$ \! $\rho_3$-pairs whose switchings realize the attachment of the 3-handles. 

More generally, we can state the following result:

\begin{proposition}\label{chiudere} Let $M$ be a compact 4-manifold with $\partial M\cong\#_r(\mathbb S^1\times\mathbb S^2)$ and let $\bar M$ be its associated closed 4-manifold. If $\G\in G_s^{(4)}$ is a gem of $M$ such that $\G_{\hat 4}$ contains $r\geq 1$ $\rho_3$-pairs of color $i\in\Delta_3$, whose switchings yield a gem of $\mathbb S^3,$ then:
\begin{itemize}
\item [(i)] if all $\rho_3$-pairs of $\G_{\hat 4}$ are $\rho_3$-pairs in $\G$, too, then their switchings yield a gem of $\bar M;$
\item [(ii)] if $r'\leq r$ of them are $\rho_4$-pairs in $\G$, then by switching all pairs a gem  of a closed 4-manifold $N$ is obtained such that:
$\bar M\cong\#_{r'}(\mathbb S^1\times\mathbb S^3)\# N.$
\end{itemize}
\end{proposition} 

\dimo
 With regard to statement (i), let us first suppose $r=1$ and let $\Gamma^\prime$ be the 5-colored graph obtained by switching the $\rho_3$-pair of color $i\in\Delta_3$ in $\G$. 
It is not difficult to see that $\Gamma^\prime$ represents a closed $4$-manifold: in fact, $\G^\prime_{\hat 4}$ represents the 3-sphere by hypothesis, while the only ${\hat c}$-residues that have been modified ($c\in\Delta_3$) still represent $\mathbb S^3$ since they are obtained by the corresponding ${\hat c}$-residues of $\G$ by switching a $\rho_2$-pair (see Section \ref{preliminaries.gems}).

\begin{figure} [h!]
   \centering
   \includegraphics[width=1\linewidth]{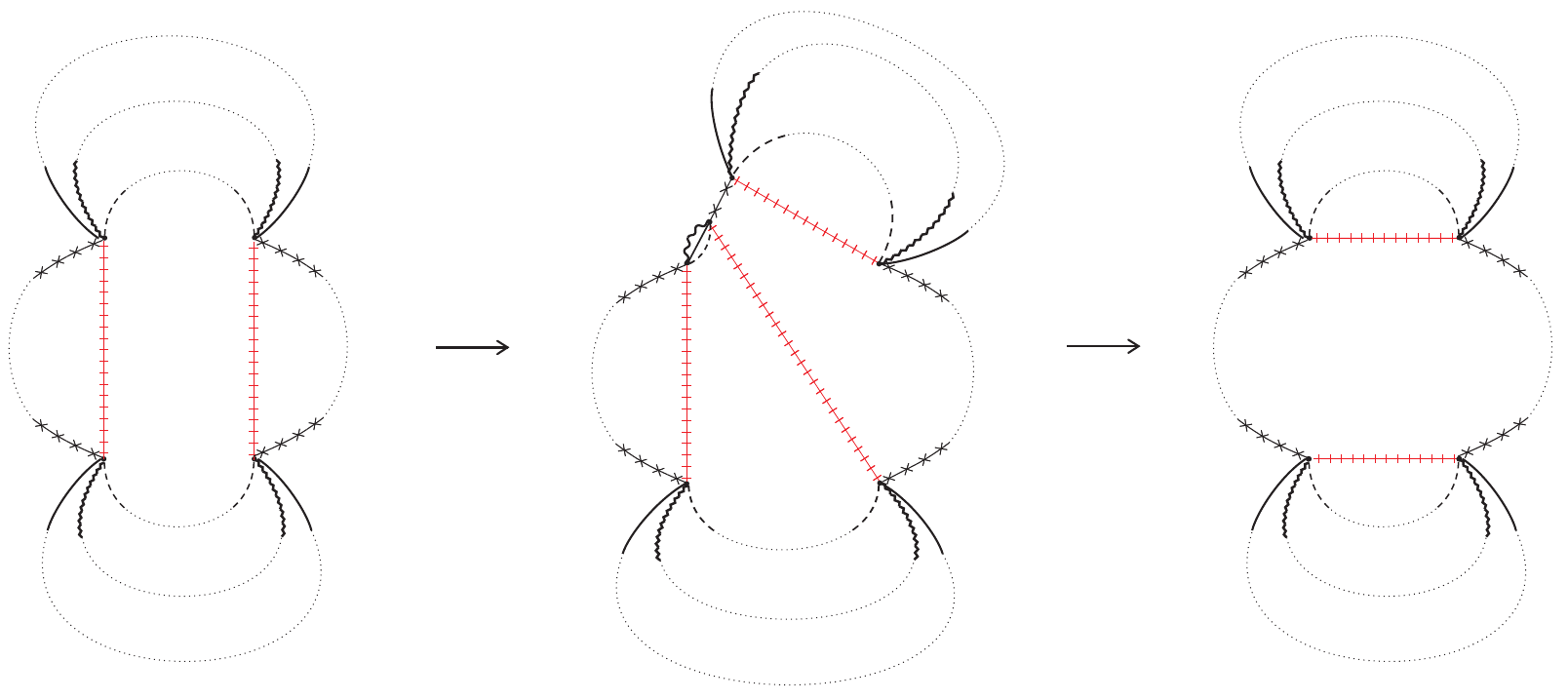}
\caption{Factorization of the switching of a $\rho_3$-pair}
    \label{fig:ro3pair}
\end{figure}

Moreover, as shown in Figure \ref{fig:ro3pair}, the switching of the $\rho_3$-pair in $\G$ is equivalent to the insertion of a $4$-colored edge on one of the edges of the pair, followed by the cancellation of a $3$-dipole. 
It is not difficult to check that the insertion of the $4$-colored edge has the topological effect of attaching a polyhedron homeomorphic to $\mathbb D^1\times\mathbb D^3$ to the boundary of $M$, without affecting its interior: in fact, it corresponds to  ``breaking'' a tetrahedral boundary face and inserting two $4$-simplices sharing the $3$-dimensional face opposite to the $4$-labelled vertex, so as to change $\partial M$ into a $3$-sphere\footnote{The same argument is used in the proof of \cite[Proposition 4.3]{Casali-Cristofori trisection bis} and also in \cite[Proposition 11(ii)]{Casali-Cristofori Kirby-diagrams} for the particular case of gems arising from Kirby diagrams.}.
As a consequence, since the subsequent cancellation of the 3-dipole does not change the represented manifold (see Section \ref{preliminaries.gems}), the closed 4-manifold $|K(\G^\prime)|$ is obtained from $M$ by attaching a 3-handle (and a 4-handle to cap off the resulting spherical boundary); by \cite{Laudenbach-Poenaru} $\G^\prime$ is thus proved to represent $\bar M.$

The case $r  > 1$ is proved in a similar way by induction, with the only difference that at each switching a $\mathbb S^1\times\mathbb S^2$ summand is ``subtracted'' from the boundary (which is represented by the $\hat 4$-residue).
\medskip

Let us now consider statement (ii) and start with the case $r=r'=1:$  let $(e,f)$ be the $\rho_4$-pair of color $i\in\Delta_3$ in $\G$ satisfying the hypotheses and $\G^\prime$ the 5-colored graph obtained by switching $(e,f).$

In order to recognize $|K(\G^\prime)|$, we factorize the switching of $(e,f)$ by some other combinatorial moves whose topological effects are known. 

First, let us insert in $\G$ a 3-dipole, which does not change the represented manifold but has the effect of disconnecting the $\{i,4\}$-colored cycle $C$ containing $(e,f)$. More precisely, let us delete an $i$-colored edge  $e'$ and a 4-colored  edge $f'$ of $C$, suitably chosen such that $e$ and $f$ belong to different connected components of $C\setminus\{e',f'\}$; then let us add two vertices, $x$ and $y$, joined by three edges of the remaining colors, and connect $x$ and $y$ with the endpoints of $e'$ (resp. $f'$) by $i$-colored (resp. 4-colored) edges, taking care to respect the bipartition of the vertices. While the obtained 5-colored graph $\tilde\G$ is still a gem of $M$, $(e,f)$ is now a $\rho_3$-pair in $\tilde\G$ and it is not difficult to prove that it satisfies the hypotheses of the proposition\footnote{Note that $\tilde\G_{\hat 4}$ is obtained from $\G_{\hat 4}$ by the addition of the 3-dipole  of vertices $x$ and $y$.}. 
Therefore, by statement (i), the switching of $(e,f)$ in $\tilde\G$ produces a new 5-colored graph $\tilde\G^\prime$ representing the closed manifold $\bar M.$  

In order to recover $\G^\prime$, it is now necessary to delete the subgraph formed by $x,y$ and paste the resulting hanging edges. However, this subgraph is a so-called {\it combinatorial handle} in $\tilde\G^\prime$, since $x$ and $y$ share the same $\{i,4\}$-colored cycle. As a consequence of  \cite[Theorem 7]{Gagliardi-Volzone}, $\G^\prime$ turns out to be a gem of a closed $4$-manifold $N$ such that $\bar M\cong |K(\tilde\G^\prime)|\cong (\mathbb S^1\times\mathbb S^3)\# N.$

More generally, if $1\leq r' \le r$, we insert $r'$ 3-dipoles in $\G$ in order to obtain a new gem $\tilde\G$ of $M$ containing $r$ $\rho_3$-pairs, then by switching them and deleting the $r'$ resulting combinatorial handles, a 5-colored graph is obtained which, by the above arguments, is proved to represent a closed $4$-manifold $N$ such that, as claimed, $\bar M\cong\#_{r'}(\mathbb S^1\times\mathbb S^3)\# N.$\qed

\bigskip

An example of the situations described in Proposition \ref{chiudere} (with $r=r^\prime=1$) can be seen in Figures \ref{fig:M4(K0,0)} (statement (i)) and \ref{fig:M4(K0,dotted)} (statement (ii)).
More precisely, the left side of Figure \ref{fig:M4(K0,0)} shows a gem of $\mathbb S^2 \times \mathbb D^2$, obtained by applying the algorithm to the 0-framed trivial knot, with the highlighted $\rho_3$-pair, whose switching gives rise to the gem of the associated closed manifold $\mathbb S^4$ depicted on the right. 
A gem of $\mathbb S^1\times\mathbb D^3$, obtained from a dotted trivial knot, can be seen on the left side of Figure \ref{fig:M4(K0,dotted)}, with a highlighted $\rho_4$-pair whose switching produces the gem on the right representing $N (=\mathbb S^4)$, so that the associated closed 4-manifold is $(\mathbb S^1\times\mathbb S^3)\# N$.  

\begin{figure} [h!]
    \centering
    \includegraphics[width=0.45\linewidth]{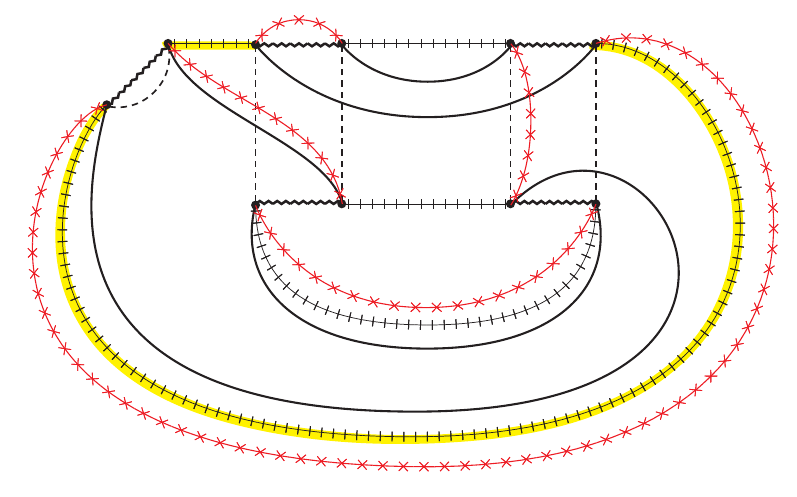}\quad\includegraphics[width=0.45\linewidth]{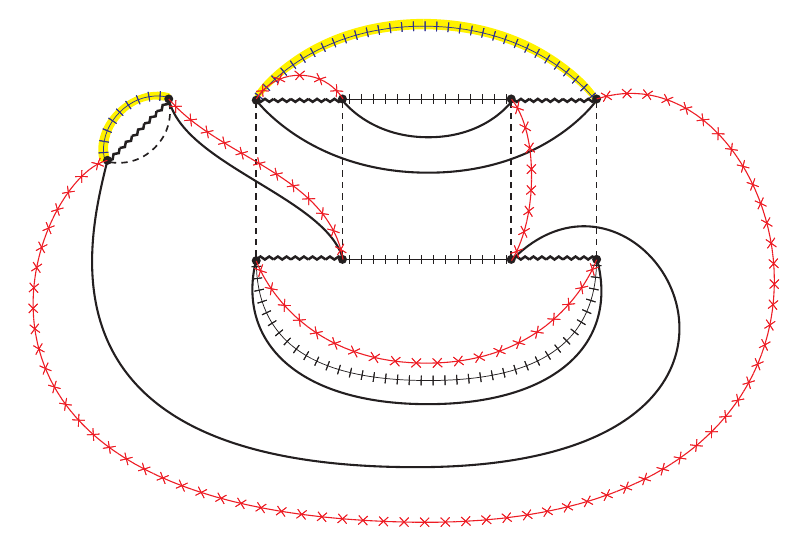}
\caption{gems of $M = \mathbb S^2 \times \mathbb D^2$ (left) and $\bar M = \mathbb S^4$ (right), arising from the $0$-framed trivial knot} 
  \label{fig:M4(K0,0)}
\end{figure}

\begin{figure} [h!]
    \centering
    \includegraphics[width=0.45\linewidth]{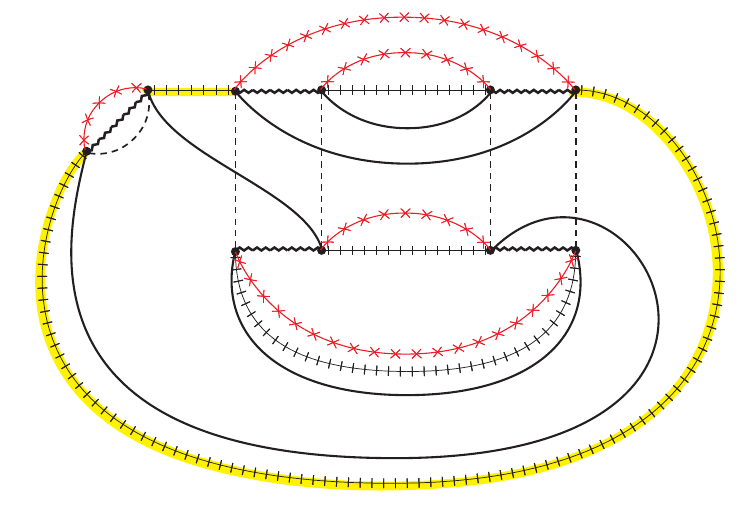}\quad \includegraphics[width=0.45\linewidth]{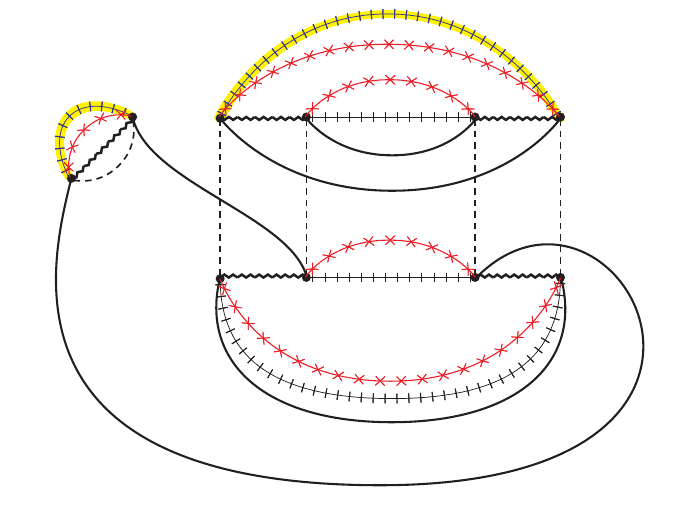}
\caption{gems of $M= \mathbb S^1 \times \mathbb D^3$ (left) and  $N(=\mathbb S^4)$ (right), so that  $\bar M = (\mathbb S^1 \times \mathbb S^3) \# N$, arising from the dotted trivial knot}
    \label{fig:M4(K0,dotted)}
\end{figure}

\bigskip

We conclude the Section by pointing out that useful applications of the described algorithm have been developed.    
In fact, it has been implemented by R.A. Burke in a software utility of {\it Regina} software package (\cite{regina}), which gives the additional advantage of allowing the manipulation of the resulting triangulations through an expressively thought-of heuristic for the 4-dimensional case.

As a consequence, in \cite{Burke-Burton-Spreer}, where a census of $4$-manifold triangulations is studied, the algorithm is applied to get triangulations of the ``canonical" PL-type of each topological type of the census, in order to perform automatically a comparison among the obtained triangulations. 

\smallskip 
Another interesting application regards the possibility of obtaining colored graphs - and hence colored triangulations -  representing exotic pairs of compact PL 4-manifolds, i.e. pairs of 4-manifolds which are TOP-homeomorphic but not PL-homeomorphic. 
(Singular) triangulations of some exotic pairs are shown in both \cite{Casali-Cristofori Kirby-diagrams} and \cite{Burke}; moreover, the first examples of triangulations of {\it corks} have been obtained in \cite{Burke}, where an initial analysis of the resulting triangulations is also presented, together with some interesting observations about particular configurations that seem to recur within.


\subsection{From gems to trisections} \label{ss.from_gems_to_trisections}
\par \noindent

In the last years, the notion of trisection has been extensively studied, also in its relationships with other manifold representation tools, such as handle-decom- positions, Lefschetz ﬁbrations, etc. In particular, in 2018, relying on the coincidence between DIFF and PL categories in dimension 4, Bell, Hass, Rubinstein and Tillmann (\cite{Bell-et-al}) faced the study of trisections via (singular) triangulations.

They make use of a so-called {\it tricoloring}, i.e. a labelling of the vertices of each 4-simplex by three colors: in fact, the pull-back of the cubical structure of the standard 2-simplex, via the natural map induced by a tricoloring, gives rise to a subdivision of the $4$-manifold into three 4-dimensional pieces.  Under suitable conditions (always realized through Pachner moves, although at the cost of increasing the number of 4-simplices) the above subdivision turns out to be a trisection.

This approach - which moves trisections from SMOOTH to PL setting - brings all the advantages of a combinatorial description, and allows to algorithmically construct trisections and estimate their genus, starting from singular triangulations.

\medskip 

Spreer and Tillmann (\cite{Spreer-Tillmann(Exp)})  performed the ﬁrst connection between trisections and gems, by applying the idea of Bell, Hass, Rubinstein and Tillmann to a special type of triangulations of simply-connected $4$-manifolds, which contain exactly 5 vertices and one 1-simplex between each pair of vertices; their dual colored graphs are called {\it simple crystallizations}.

In this case, an obvious tricoloring of the vertices can be considered. 
Spreer and Tillmann proved that, for at least one simple crystallization of each $4$-manifold among $\mathbb{CP}^2$, $\mathbb{S}^2 \times \mathbb{S}^2$ and the $K3$-surface, the associated subdivision is actually a trisection.
 
 Moreover, they computed the genus of the central surface and pointed out that it had to be the minimum one, since in the case of simple crystallizations all 4-dimensional pieces of the subdivision are disks.   

As a consequence, they completed the computation of the trisection genus for all the so-called ``standard" simply-connected 4-manifolds, i.e. for the connected sums of $\mathbb{CP}^2$ - possibly with reversed orientation -, $\mathbb{S}^2 \times \mathbb{S}^2$ and the $K3$-surface; and for all such manifolds the trisection genus turns out to coincide with the second Betti number.

\bigskip

In this Section we will see how Spreer and Tillmann’s construction can be generalized, so as to yield  trisections (directly or indirectly) from \underline{each} colored triangulation of a closed 4-manifold, both in the orientable and non-orientable case.  
By the way, also compact 4-manifolds with connected boundary will be taken into consideration, by introducing a suitable extension of the notion of trisection to the boundary case (diﬀerent from the one studied by Castro, Gay and Pinzon in \cite{Castro-Gay-Pinzon}).

Note that, in order to study exotic structures, it is of fundamental importance to go beyond the class of simple crystallizations also in the simply-connected setting: in fact it is not diﬃcult to prove the existence of inﬁnitely many simply-connected PL 4-manifolds, TOP-homeomorphic to the same standard simply-connected PL 4-manifold,  and not admitting simple crystallizations (see \cite[Prop. 32]{Casali-Cristofori ElecJComb 2015}).


\subsubsection{Gem-induced trisections}  \label{sss Gem-induced trisections} \par \noindent

In  \cite{Casali-Cristofori gem-induced} (resp. \cite{Casali-Cristofori trisection bis}), the ideas of \cite{Bell-et-al} are applied to all colored triangulations of compact orientable (resp. non-orientable) 4-manifolds with empty or connected boundary, so as to generalize Spreer and Tillmann’s construction. 

As a starting point, it is proved that each gem of a compact $4$-manifold $M$ with empty or connected boundary, both in the orientable and non-orientable case, induces a decomposition of $M$ which ``resembles" a trisection:  there are three 4-dimensional ``pieces", two of which are handlebodies, while the third one is a collar of the boundary (or a 4-disk, if the boundary is empty); moreover, all these sub-manifolds
intersect into a closed surface, but only two of the pairwise intersections are ensured to be 3-dimensional handlebodies. 

\begin{theorem}  \ {\rm (\cite{Casali-Cristofori gem-induced}, \cite{Casali-Cristofori trisection bis})} \  \label{th: gem-induced trisection}
For each 5-colored graph $\Gamma \in G_s^{(4)}$ representing a compact $4$-manifold with empty (resp. connected boundary) $M$ and for 
each cyclic permutation $\varepsilon= (\e_0,\e_1,\e_2,\e_3, \e_4=4)$ of $\Delta_4,$ 
a  triple $\mathcal  T(\Gamma, \varepsilon) =(H_{0},H_{1},H_{2})$ of submanifolds of $M$ is constructed, such that
 \begin{itemize}
 \item [(a)]  $M = H_{0}\cup H_{1}\cup H_{2}$ and $H_{0}, H_{1}, H_{2}$ have pairwise disjoint interiors 
 \item [(b)]  $H_{1},H_{2}$ are $4$-dimensional handlebodies; $H_{0}$ is a $4$-disk \  (resp. is homeomorphic to $\partial M \times [0,1]$)
\item [(c)] $H_{01}=H_{0}\cap H_{1},$  $H_{02}=H_{0}\cap H_{2}$ are $3$-dimensional handlebodies
 \item [(d)] $\Sigma = H_{0}\cap H_{1}\cap H_{2}$  is a closed connected surface. 
 \end{itemize}
 
 Moreover, if $H_{12}=H_1\cap H_2$ is a 3-dimensional handlebody, too, then all the above handlebodies, as well as the surface $\Sigma$, are orientable or not according to the orientability of $M$; in the first case $\Sigma$ has genus $\rho_{\varepsilon}(\G_{\hat 4})$, while in the second one it has genus $2\rho_{\varepsilon}(\G_{\hat 4}).$ 
\end{theorem}

{\it Sketch of the proof:}
\par \noindent 
According to \cite{Bell-et-al}, for each tricoloring on the vertices of a triangulation, the possible trisecting submanifolds are obtained as pull-backs of the cubical structure of the standard 2-simplex, via the natural map induced by the tricoloring itself. 

If the triangulation associated to a colored graph is considered, where all 4-simplices have one vertex of each color from $0$ to $4$, a tricoloring is easily obtained by giving color red to both the $\e_0$- and $\e_2$-labelled vertices, color green to both the $\e_1$- and $\e_3$-labelled vertices, and color blue to the vertex labelled by $\e_4=4$ (which can be actually assumed to be unique, since the boundary is either empty\footnote{Obviously, in the closed case, the role of color $4$ can be played by any color $c$ such that $\Gamma_{\hat c}$ is connected. In this paper, we always use $4$ for simplicity, both in the closed and boundary case.} or connected: see Theorem \ref{Theorem_gem}). 

Hence, one of the 4-dimensional pieces corresponds to the (truncated\footnote{Only in the boundary case, a small open neigbourhood of the singular ($4$-labelled) vertex has to be removed, according to gem theory.}) cone over the link of the unique blue vertex in the first barycentric subdivision of the triangulation (which is a 3-sphere in the closed case, and the boundary 3-manifold  otherwise), while the other two are 4-dimensional handlebodies, obtained as regular neighbourhoods of the 1-dimensional subcomplexes of $K(\Gamma)$ generated by its red and green vertices, respectively. 

As regards the pairwise intersections, which are 3-dimensional submanifolds, it is not difficult to prove that two of them meet any 4-simplex in triangular prisms, while the third (corresponding to the blue singleton) meets it in a cube. Moreover, any two of them meet in a unique square, and this square is exactly the intersection between the 4-simplex and the central surface.

\begin{figure}   [h!]
     \centering
    \includegraphics[width=1\linewidth]{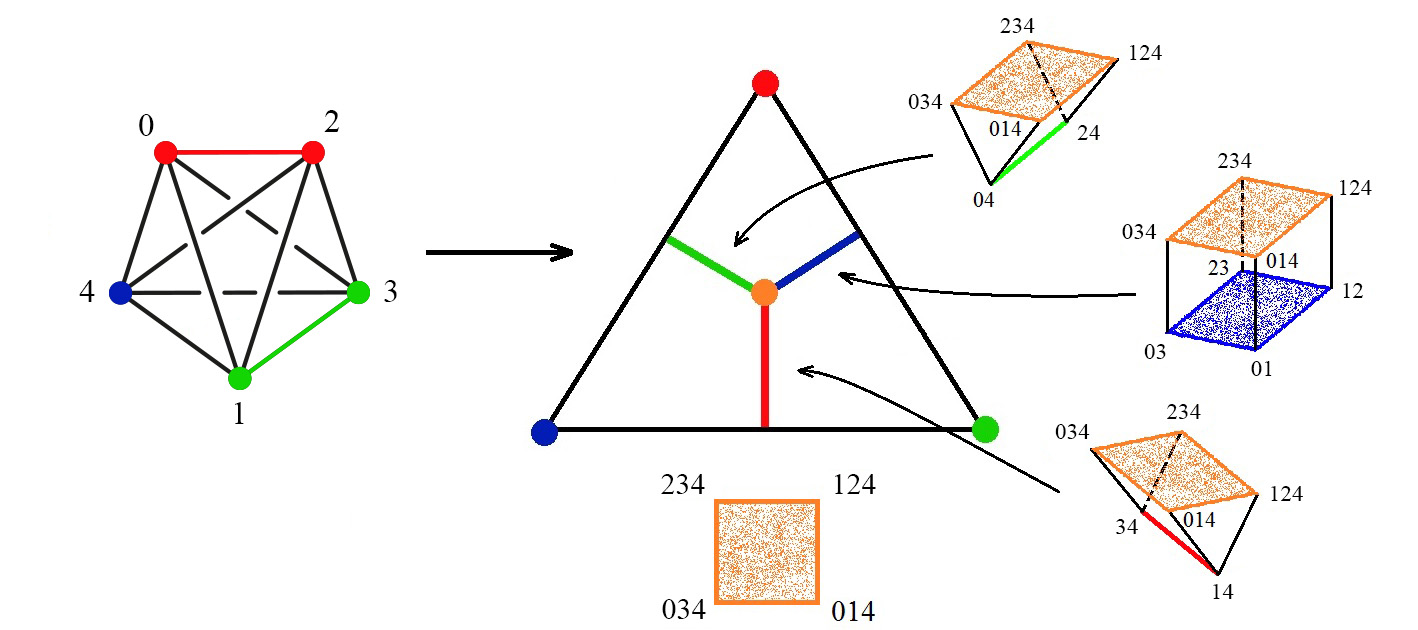}
    \caption{The intersections of the 3-dimensional pieces and the central surface of $\mathcal  T(\Gamma, \varepsilon)$ with any 4-simplex}
     \label{fig:trisec_pieces}
 \end{figure}

See Figure \ref{fig:trisec_pieces} (redrawn from \cite{Spreer-Tillmann(Exp)}), where $\e=(0,1,2,3,4)$ is assumed:
the vertices of the building blocks of the trisecting submanifolds, which are barycenters of faces, are indicated by the labels of the spanning vertices. 

\bigskip
It is not difficult to check that both the triangular prisms corresponding to color red and those corresponding to color green always collapse to 1-dimensional subcomplexes, and hence are 3-dimensional handlebodies.

Moreover, the triple intersection $\Sigma= H_0\cap H_1\cap H_2$ is a closed connected surface, whose Euler characteristic can be easily computed. \qed

\bigskip
Obviously, the construction is interesting when also the third 3-dimensional intersection turns out to be a handlebody, giving rise to an effective generalization to the connected boundary case of the notion of trisection, induced by graphs encoding manifolds: 

\begin{definition}   \ {\rm (\cite{Casali-Cristofori gem-induced}, \cite{Casali-Cristofori trisection bis})} \  
{\em Let $M$ be a compact $4$-manifold with empty or connected boundary. 
A {\it gem-induced trisection} of $M$ is a decomposition 
$\mathcal  T(\Gamma, \varepsilon) =(H_{0},H_{1},H_{2})$ such that $H_{12}=H_1\cap H_2$ is a 3-dimensional handlebody, $\Gamma \in G_s^{(4)}$ being a 5-colored graph representing $M$ and $\e=(\e_0,\e_1,\e_2,\e_3,4)$ a cyclic permutation of $\Delta_4$.}
\end{definition} 

\begin{remark}
{\em Note that, in the closed case, gem-induced trisections belong to a very significant type of trisections, where one of the pieces is a $4$-disk. Such trisections have been widely studied, and are conjectured to exist for all closed simply-connected $4$-manifolds: see, for example, \cite{Meier}, \cite{Lambert-Cole-Meier} and \cite{Meier-Schirmer-Zupan}. 

On the other hand, when the boundary is non-empty, the obtained decomposition of $M$ is different  from that introduced in  \cite{Castro-Gay-Pinzon}, which is commonly intended as a trisection in the boundary case: in fact, it doesn't involve compression bodies and open book decompositions, but only handlebodies and Heegaard splittings, according to a suggestion by Rubinstein and Tillmann in  \cite[Section 4.5]{Rubinstein-Tillmann}. }
\end{remark}

\bigskip
From the combinatorial point of view, a condition ensuring that a given gem admits a gem-induced trisection has been determined: 

\begin{proposition} \label{CS gem-induced trisections} \ {\rm (\cite{Casali-Cristofori gem-induced}, \cite{Casali-Cristofori trisection bis})} \ 
Let $M$ be a compact $4$-manifold with empty or connected boundary and $\G\in G_s^{(4)}$ a gem of $M$ of order $2p$;  if  there exists an ordering $(e_1,\ldots,e_p)$ of the 4-colored edges of $\G$ such that for each $j\in\{1,\ldots,p\}$:
$$\begin{aligned} \text{there exists } i\in\Delta_3 \text{ such that} & \ \text{all 4-colored edges of  the }  \\ \{4,i\}\text{-colored cycle} \ \text{containing } e_j  & \ \text{belong to the set }\{e_1,\ldots,e_j\},\end{aligned}\hskip 15pt (*)$$     
then $\mathcal T(\G, \varepsilon)$ is a gem-induced trisection of $M$, for each cyclic permutation $\e$ of $\Delta_4.$       
\end{proposition} 

\medskip

{\it Sketch of the proof:}
\par \noindent 
It is easy to prove that the union of cubes, constituting the third 3-dimensional intersection $H_{12}=H_1\cap H_2$ of $\mathcal T(\G, \varepsilon)$, always collapses to the $2$-dimensional complex $Q(\G,\varepsilon)$ consisting of the union of their bottom faces.  
On the other hand, each such face corresponds to a 4-colored edge of $\Gamma$, and its edges are dual to bicolored cycles of $\Gamma$ involving color 4. 

\begin{figure}[h!]
\centering
\scalebox{0.8}{\includegraphics{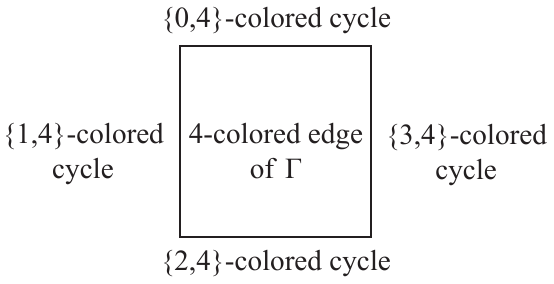}}
\caption{the square, corresponding to a $4$-colored edge of $\G$, constituting $Q(\G,\varepsilon)$}
\label{figHsquare}
\end{figure}

In order to prove that, under the hypothesis of the statement, $Q(\G,\varepsilon)$ further collapses to a graph, first note that condition (*), for $j=1$, means that $e_1$ belongs,  for a suitable color $i$, to a $\{4,i\}$-colored cycle of $\G$ of length two;  hence, the square  in $Q(\G,\varepsilon)$ corresponding  to $e_1$ has a free edge from which it can be collapsed (see Figure \ref{figHsquare}, where $\e=(0,1,2,3,4)$ is assumed).

Now, it is not difficult to check that an ordering of all $4$-colored edges of $\G$  satisfying condition (*) corresponds exactly to a sequence of elementary collapses of all squares of $Q(\G,\varepsilon)$.\qed 

\bigskip

It is worthwhile to note that gem-induced trisections fit Gay and Kirby's definition of trisection - hence allowing a direct estimation of the trisection genus - only for a suitable class of closed orientable $4$-manifolds, characterized by admitting a handle decomposition with no 3-handles:  

\begin{proposition}\label{trisection_vs_gem-induced}  \  {\rm \cite{Casali-Cristofori trisection bis} } \ 
Let $M$ be a compact $4$-manifold with empty or connected boundary. 
A gem-induced trisection of $M$ is a trisection if and only if  $M$ is closed and orientable. 

\noindent Moreover, a closed orientable $4$-manifold admits a gem-induced trisection if and only if it admits a handle decomposition lacking in 3-handles (or in 1-handles).
\end{proposition}

The positive aspect of the above result is that the involved class possibly comprehends all simply-connected closed 4-manifolds, according to the famous Kirby problem n. 50 (see \cite[Section 7.3]{[M]}).

On the other hand, Proposition \ref{trisection_vs_gem-induced} shows that it is impossible to face the whole class of closed $4$-manifolds, looking for an estimation of their trisection genus, directly via gem-induced trisections.  

Nevertheless, starting from gem-induced trisections,  
two ``indirect" approaches to trisections of closed (orientable and non-orientable) 4-manifolds can be applied, starting from gem-induced trisections. They will be the object of the following Sections.


\subsubsection{First indirect approach} \label{ss: First indirect approach} \par \noindent

The first indirect approach to trisections of closed 4-manifolds via gem theory relies on the particular type of extension to the boundary case performed by gem-induced trisections, where one of the $4$-dimensional pieces is a collar of the boundary.   In fact, if the boundary of a compact $4$-manifold  $M$ is homeomorphic to a connected sum of sphere bundles over $\mathbb S^1$ (both in the orientable and non-orientable case), it is not difficult to prove that any  gem-induced trisection of $M$ naturally gives rise to a trisection of the associated closed 4-manifold $\bar M$, i.e. the one uniquely obtained by gluing a $4$-dimensional handlebody along $\partial M$, according to the already cited theorem by Laudenbach and Poenaru (\cite{Laudenbach-Poenaru}), together with its non-orientable version (\cite{Miller-Naylor}).  

Since the intersecting surface is not affected,  an upper bound for the trisection genus of the closed 4-manifold $\bar M$ is obtained, through the so-called \emph{G-trisection genus} of $M$ (denoted by $g_{GT} (M)$), that is the minimum genus of a central surface in a gem-induced trisection: 

\begin{theorem}\label{trisection_from_gem-induced} \ {\rm (\cite{Casali-Cristofori trisection bis})} \ 
\ Let \ $M$ \ be \ a \ compact $4$-manifold \ with \ boundary  \  \  $\partial M \cong \#_m(\mathbb S^1 \otimes \mathbb S^2)$,  with $m>0$.           
If \ $M$\  admits\  a\  gem-induced\  trisection,\  then\  $\bar M \cong M \cup {\mathbb Y}_m^{(\sim)}\!$   admits a trisection with the same central surface.  

As a consequence, $$g_T(\bar M) \le g_{GT} (M).$$
\end{theorem}

\medskip 

\begin{definition}
{\em A trisection of a closed (orientable or non-orientable)  4-manifold $\bar M$ is said to {\it arise from a colored triangulation} if it is either a gem-induced trisection of $\bar M$ or is obtained from a gem-induced trisection of $M$ (such that $\bar M \cong M \cup {\mathbb Y}_m^{(\sim)}$) according to Theorem \ref{trisection_from_gem-induced}.     }
\end{definition} 

\bigskip

In \cite{Casali-Cristofori gem-induced} and \cite{Casali-Cristofori trisection bis}), various examples of closed $4$-manifolds are shown, for which the estimation of the trisection genus via gems (either directly, through G-trisection genus, or indirectly through the trisection genus of an associated bounded manifold) is sharp. 
Moreover, since this type of estimation is proved to be subadditive with respect to connected sum,  
trisections arising from colored triangulations turn out to minimize the trisection genus for a wide class of (orientable and non-orientable) $4$-manifolds: 
 
\begin{proposition} \label{calculations} \ {\rm (\cite{Casali-Cristofori trisection bis})} \  Let \ 
$\bar M\, \cong_{PL}\,(\#_p\mathbb {CP}^2)\,\#\,(\#_{p^{\prime}}(-\mathbb {CP}^2))\, \#\, (\#_q(\mathbb S^2\times \mathbb S^2)) \, \#  \, (\#_r K3) \,\#\,(\#_s(\mathbb S^1 \otimes \mathbb S^3))\, \# (\#_t\mathbb R \mathbb P^4)\,\# \, (\#_u  (\mathbb S^2 \times \mathbb {RP}^2 )),$  with $p,p^{\prime},q, r,s,t,u \geq 0$.  
Then, its trisection genus \ $ g_T(\bar M)  = (p + p^{\prime} + 2q + 22r)+s+2t+3u$ \ 
is realized by a trisection arising from a colored triangulation.
\end{proposition}

\medskip

\subsubsection{Second indirect approach}  \par \noindent

The second ``indirect" approach to trisections via colored triangulations is used when the subdivision of the $4$-manifold induced by a given gem has one of the 3-dimensional pieces (namely $H_{12}$, according to the notations of Section \ref{sss Gem-induced trisections}) that does not collapse to a graph.

\medskip 
In a recent paper by Martini and Toriumi (\cite{Martini-Toriumi}), the idea of performing {\it stabilizations} along all edges of a fixed color is introduced, in order to obtain a trisection from each gem of a closed $4$-manifold.   

\begin{figure}   [h!]
     \centering
     \includegraphics[width=0.5\linewidth]{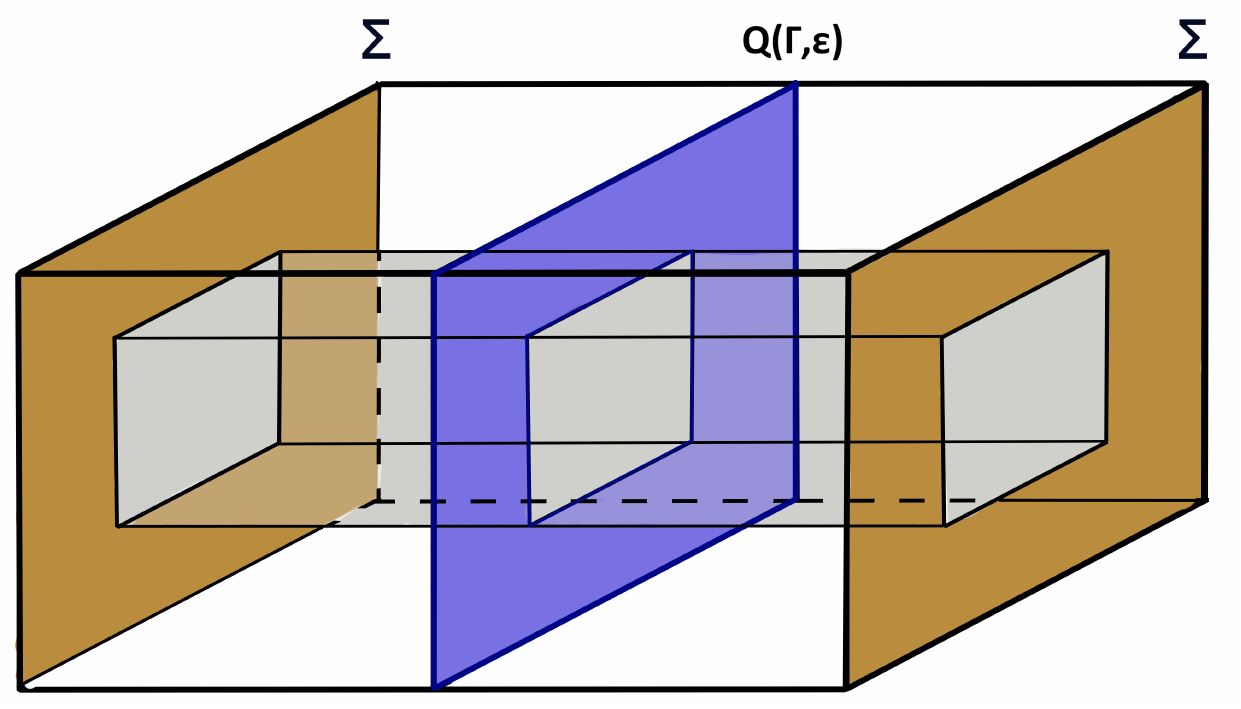}
     \caption{stabilization along a $4$-colored edge}
     \label{fig:carving_edge}
 \end{figure}

By translating their notations into those of Section \ref{sss Gem-induced trisections}, what they prove is that, by carving a regular neighborhood of a $4$-colored edge $\bar e$,\footnote{Note that, if $\Gamma_{\hat 4}$ is connected - as it is assumed in Section \ref{sss Gem-induced trisections} - the case of $4$-colored edges connecting distinct $\hat 4$-residues never occurs. On the contrary, Martini and Toriumi also take into account this situation, where carving along an edge simply corresponds to a connected sum of ``pieces" of the decomposition.} 
a $1$-handle is added to $H_0$, without affecting the $4$-dimensional handlebodies $H_1$ and $H_2$,  while the genus of the intersecting surface $\Sigma$ increases by one; with regard to the $3$-dimensional pieces, the operation consists in adding a $1$-handle to $H_{01}$ and $H_{02}$, while adding a $2$-handle to $H_{12}$.  
Figure \ref{fig:carving_edge} shows the effect of stabilization on the cubes that are building blocks of $H_{12}$, specifically on the pair of adjacent cubes corresponding to the 4-colored edge $\bar e$. In particular, in the middle square, which is part of the $2$-dimensional spine $Q(\G,\varepsilon)$ of $H_{12}$, a “hole” is created, which ensures the collapse of the square to its boundary.

By performing the operation on every $4$-colored edge of a gem $\Gamma \in G_s^{(4)}$, the collapse of $Q(\G,\varepsilon)$ to a $1$-dimensional subcomplex is guaranteed. In this way, Martini and Toriumi prove the existence of trisections arising from $\Gamma$, at the cost of increasing the genus of the central surface (since the final genus is the sum of the regular genus of $\Gamma_{\hat 4}$ and the number of stabilizations).

\medskip

In \cite{Casali-Cristofori trisection ter}, the idea of Martini and Toriumi is applied to gems (belonging to  $G_s^{(4)}$) representing compact $4$-manifolds with empty or connected boundary, with the purpose of minimizing the number of stabilizations required to yield a trisection.

Both in the closed case and in the case of boundary homeomorphic to a connected sum of sphere bundles over $\mathbb S^1$ (when the first indirect approach to trisections can be applied), the  genus of the obtained trisection turns out to be less or equal to the regular genus of the starting gem: 

\begin{theorem}\label{number_stabilizations} {\rm (\cite{Casali-Cristofori trisection ter})} \ 
Let $\G\in G_s^{(4)}$ be a gem of a compact $4$-manifold $M$ with $\partial M \cong \#_m(\mathbb S^1 \otimes \mathbb S^2)$,  $m \ge 0$.   
For each cyclic permutation $\varepsilon$ of $\Delta_4$, the  triple $\mathcal  T(\Gamma, \varepsilon)$ gives rise to a trisection of $\bar M \cong M \cup{\mathbb Y}_m^{(\sim)}$  by applying $k$  stabilizations, with $0 \le k \le \rho_\e(\Gamma) - \rho_{\e}(\Gamma_{\hat 4})$.
Hence: 
$$g_T(\bar M) \le \rho_\e(\Gamma)$$
\end{theorem}

{\it Sketch of the proof:}
\par \noindent 
First, let us consider the (possibly disconnected) $3$-colored graph $\G_{\{\e_0,\e_{3},\e_4=4\}}$, which regularly embeds in $\mathbb S^2$ (as happens for each 3-residue of a gem in dimension $ \ge 3$).  
This fact allows ordering all $4$-colored edges of $\Gamma$, so that, after performing a stabilization along the first $g_{\e_0 \e_3}-g_{\e_0 \e_3 4}$ edges, the remaining ones satisfy condition (*) of Proposition \ref{CS gem-induced trisections}. 

\par Further, well-known combinatorial formulas in gem theory (see for example \cite[Proposition 3.1]{generalized-genus}) enable to check that 
$g_{\e_0 \e_3}-g_{\e_0 \e_3 4} = \rho_\varepsilon(\Gamma) - \rho_\varepsilon(\Gamma_{\hat 4})$.

Hence, a gem-induced trisection of $M$ is obtained, with central surface of genus:  $\rho_\varepsilon(\Gamma_{\hat 4}) + (\rho_\varepsilon(\Gamma) - \rho_\varepsilon(\Gamma_{\hat 4}))$. 

Finally, the thesis concerning $\bar M$ follows, by possible application of Theorem \ref{trisection_from_gem-induced} (in  case $m>0$). \qed 

\medskip

As a direct consequence, the following statement holds:

\begin{theorem}\label{trisection_vs_regular-genus} \ {\rm (\cite{Casali-Cristofori trisection ter})} \ 
For each compact $4$-manifold $M$ with $\partial M \cong \#_m(\mathbb S^1 \otimes \mathbb S^2)$, $m \ge 0$,    then: 
$$g_T(\bar M) \le \mathcal G(M)$$
where $\bar M \cong M \cup {\mathbb Y}_m^{(\sim)}$.  
\end{theorem}

\medskip

In particular, in the closed case, the regular genus turns out to be an upper bound for the trisection genus: 

\begin{theorem}\label{trisection_vs_regular-genus(closed case)} \ {\rm (\cite{Casali-Cristofori trisection ter})} \ 
For each closed $4$-manifold $M$,  then: 
$$g_T(M) \le \mathcal G(M).$$
\end{theorem}

\medskip

\subsection{From Kirby diagrams to trisections via gems} 
\label{ss.from_Kirby_to_trisections}
\par \noindent

In the present Section, we will show how the results and constructions in Section \ref{ss.from_gems_to_trisections} take advantage of those of Section \ref{ss.from_Kirby_to_gems}, as regards the computation of the trisection genus of each closed orientable 4-manifold. 

\medskip

The starting point is the fact, proved in \cite{Casali-Cristofori gem-induced}, that all 5-colored graphs obtained through the algorithmic construction described in Section \ref{ss.from_Kirby_to_gems}, starting from Kirby diagrams (see Theorem \ref{t.Kirby-diagram}), as well as the simplified ones arising from framed links (see Theorem \ref{regular-genus&gem-complexity(framed)}), admit gem-induced trisections:  

\begin{theorem} \label{th_M4(L,c)} \ {\rm (\cite{Casali-Cristofori gem-induced})} \  Let $M$ be a compact $4$-manifold with empty or connected boundary admitting a handle decomposition with no $3$-handles.  
\smallskip

\noindent 
Then, from each (connected and with dotted components in ``good position") Kirby diagram $(L,d)$ of $M$, a gem-induced trisection of $M$ can be algorithmically constructed, whose central surface has genus $s+1$, $s$ being the crossing number of the chosen diagram.

\noindent Furthermore, if $(L,d)$ has no dotted components (equivalently, if $1$-handles are missing
), then a gem-induced trisection of $M$ can be obtained,  whose central surface has genus $m_\alpha$, where $m_{\alpha}$ is the number of $\alpha$-colored regions in a chess-board coloration of the diagram, $\alpha$ being the color of the unbounded region.
\end{theorem} 

{\it Sketch of the proof:}
\par \noindent  
The combinatorial properties of the graphs obtained through the constructions of Section \ref{ss.from_Kirby_to_gems} 
always allow ordering the set of their 4-colored edges, so as to satisfy the collapsing condition (*) of Proposition \ref{CS gem-induced trisections}. 

With regard to the graph $\Gamma(L,d)$, it is sufficient to consider the following sequence: 
\begin{itemize}
\item   
for increasing values of $j \in \{m+1, \dots, l\}$, the triad of edges between vertices $P_{2r}^{(j)}$ and $P_{2r+1}^{(j)}$, for $r=0,1,2$,  of the quadricolor corresponding to the $j$-th framed component, followed by the $4$-colored edges corresponding to the sequence  of highlighted segments of arcs on the $j$-th (framed) component, and finally by the other $4$-colored edges as they are met    
along that component; 
\item  
for each $j \in \{1, \dots, m\}$, all $4$-colored edges corresponding to the $j$-th dotted component, in any order, except those of the $\{1,4\}$-cycle containing $v_j$ and $v_j^\prime$,  then the $4$-colored edge between  $v_j$ and $v_j^\prime$, and finally  the remaining $4$-colored edge. 
\end{itemize}

As regards the 5-colored graph $\tilde \Omega(L,d)$ (whose $\hat 4$-residue is obtained from that of  
$\Gamma(L,d)$, when dotted components are missing, via generalized dipole eliminations), it is sufficient to consider, for increasing values of $j \in \{1, \dots, l\}$, first the triad of edges between vertices $P_{2r}^{(j)}$ and $P_{2r+1}^{(j)}$, for $r=0,1,2$,  
and then the other $4$-colored edges  as they are  met along the $j$-th component. 

Hence, in both cases, a gem-induced trisection exists, whose genus can be directly computed via Proposition \ref{order-regular genus Gamma}. \qed

\bigskip

As a consequence, an estimation of the trisection genus $g_T(\bar M)$ is obtained for \underline{each} closed orientable 4-manifold $\bar M$, in terms of the combinatorial properties of a Kirby diagram representing  it.   

\begin{theorem}\label{trisection_from_Kirby-diagram} \  {\rm (\cite{Casali-Cristofori trisection bis})} \ 
\begin{itemize}
\item[(i)]
For each closed orientable $4$-manifold $\bar M$,  
$$g_T(\bar M) \le s+1,$$  $s$ being the crossing number of a (connected and with dotted components in ``good position") Kirby diagram representing $\bar M$. 
\item[(ii)] Furthermore, if $\bar M$ admits a handle decomposition lacking in $1$-handles, then 
$$g_T(\bar M) \le m_\alpha,$$ where
$m_{\alpha}$ is  the number of $\alpha$-colored regions in a chess-board coloration of a (connected and with no dotted component) Kirby diagram representing  $\bar M$, $\alpha$ being the color of the unbounded region.
\end{itemize}

\end{theorem}

{\it Sketch of the proof:}
\par \noindent  
As recalled in Section \ref{ss:Kirby diagrams}, if $(L,d)$ is a Kirby diagram (with the required properties as regards connectedness and ``good position" of 1-handles) representing $\bar M$, $(L,d)$ also represents the compact 4-manifold $M$ consisting only of the union of the 0-, 1- and 2-handles of the associated handle-decomposition. 

Now, Theorem \ref{th_M4(L,c)} ensures the existence of a gem-induced trisection of $M$ whose central surface has genus $s+1$ (or, better, $m_\alpha$ if $(L,d)$ contains no dotted component). 

Hence, the thesis follows from Theorem \ref{trisection_from_gem-induced}. \qed

\medskip

In case of disconnected Kirby diagrams, instead of performing Reidemeister moves in order to apply Theorem \ref{trisection_from_Kirby-diagram}, it is also possible to take advantage of the well-known fact that the represented compact $4$-manifold is the boundary connected sum of the ones represented by each connected component.  As a consequence, by repeated application of Theorem \ref{trisection_from_Kirby-diagram} (statement (i) or (ii), according to the case), together with subadditivity of the genus of trisections arising from gems (see \cite[Proposition 3.8(iv)]{Casali-Cristofori trisection bis} and the proof of \cite[Proposition 4.4]{Casali-Cristofori trisection bis}), the following result can be stated:  

\begin{corollary}\label{disconnected_diagrams} \  {\rm (\cite{Casali-Cristofori trisection bis})} \ 
Let  $\bar M$ be a closed orientable $4$-manifold and $(L,d)$ a Kirby diagram of $\bar M$ with $c$ connected components and whose dotted components - if any - are in good position. Then:
\begin{itemize}
\item[(i)]
 $$g_T(\bar M) \le s+c ,$$  $s$ being the crossing number of $(L,d)$. 
\item[(ii)]
Furthermore, if  $(L,d)$ has no dotted components,
then 
$$g_T(\bar M) \le m_\alpha+c-1,$$
where $m_{\alpha}$ is the number of $\alpha$-colored regions in a chess-board coloration of $(L,d)$, $\alpha$  being the color of the  unbounded region.
\end{itemize}
\end{corollary}


\subsubsection{From Kirby diagrams to trisection diagrams via gems} \label{ss.from_Kirby_to_trisection_diagrams}
\par \noindent

In \cite{Martini-Toriumi}, Martini and Toriumi analyze the decomposition generated by their procedure to yield trisections from any gem, and point out that it gives rise to a collection of trisection diagrams of the appropriate genus, since they obtain a redundant amount of attaching curves of the 3-dimensional intersections: see \cite[Section 4.6]{Martini-Toriumi} for details. 

\smallskip

By adapting the idea of Martini and Toriumi to the case of a gem $\Gamma$ satisfying condition (*) of Proposition \ref{CS gem-induced trisections} 
(and therefore admitting a gem-induced trisection $\mathcal T(\Gamma, \e)$), a trisection diagram can actually be identified on the central surface of the induced trisection,  which turns out to coincide with 
the embedding surface $F_\e(\Gamma_{\hat 4})$ of the graph without color 4. 
$F_\e(\Gamma_{\hat 4})$ has genus $\rho_{\e}(\Gamma_{\hat 4})$ in the orientable case and $2\rho_{\e}(\Gamma_{\hat 4})$ in the non-orientable one, and the systems of curves drawn on it and constituting a trisection diagram are nothing but suitable bicolored cycles of the given gem.

This procedure succeeds not only in the closed case, but also in the case of boundary homeomorphic to a connected sum of sphere bundles over $\mathbb S^1$,  where the gem-induced trisection of the bounded 4-manifold $M$ uniquely identifies - according to Section \ref{ss: First indirect approach} - a trisection of the associated closed 4-manifold $\bar M$. 

\medskip

\begin{proposition} \label{trisection_diagrams_gem-induced} \ 
Let $\Gamma\in G_s^{(4)}$ be a gem of a compact $4$-manifold $M$ with $\partial M \cong \#_m(\mathbb S^1 \otimes \mathbb S^2)$, $m \ge 0$, satisfying condition (*) of Proposition \ref{CS gem-induced trisections}. 

Then, for each cyclic permutation $\varepsilon=(\e_0,\e_1, \e_2, \e_3, 4)$ of $\Delta_4$, a trisection diagram for the associated closed $4$-manifold $\bar M$ 
is given by the following three systems of curves $\{ \alpha, \beta, \gamma\}$ on 
$F_\e(\Gamma_{\hat 4}):$   
\begin{itemize}
\item[-] the $\alpha$ curves (resp. $\beta$ curves) consist of all $\{\e_0,\e_2\}$-colored cycles (resp. $\{\e_1,\e_3\}$-colored cycles) of $\Gamma$, but those corresponding to a maximal tree of the 1-dimensional subcomplex $K_{\e_1 \e_3}$ (resp. $K_{\e_0 \e_2}$), generated by the vertices labelled $\e_1$ and $\e_3$ (resp. $\e_0$ and $\e_2$) of $K(\Gamma).$    
\item[-] the $\gamma$ curves are a subset of the $\{c,4\}$-colored cycles of $\Gamma$, with $c \in \Delta_3$, which can be combinatorially identified by following - via condition (*) - the collapsing sequence of the 3-dimensional handlebody $H_{12}$ (according to the notations of Section \ref{sss Gem-induced trisections}). 
\end{itemize} 
\end{proposition} 

\dimo
According to Theorem \ref{th: gem-induced trisection}, the gem-induced trisection $\mathcal  T(\Gamma, \varepsilon)$ of 
$M$ is formed by the triple $(H_{0},H_{1},H_{2})$, where $H_{1}$ (resp. $H_{2}$) is a regular neighbourhood of $K_{\e_0 \e_2}$ (resp. $K_{\e_1 \e_3}$), while  $H_{0}$ is either the $4$-disk obtained by coning over the triangulation $K(\Gamma_{\hat 4})$ of the $3$-sphere (if $m=0$) or a collar of the triangulation $K(\Gamma_{\hat 4})$ of $\partial M \cong \#_m(\mathbb S^1 \otimes \mathbb S^2)$ (if $m>0$). 

Moreover, in case $m>0$, $\mathcal  T(\Gamma, \varepsilon)$ naturally gives rise - according to Theorem \ref{trisection_from_gem-induced} - to a trisection $\bar{\mathcal  T} = (\bar H_0, H_{1},H_{2})$ of the closed $4$-manifold $\bar M \cong M \cup {\mathbb Y}_m^{(\sim)}$, where $\bar H_0= H_{0} \cup {\mathbb Y}_m^{(\sim)}$, with the same central surface $\Sigma$ as $\mathcal  T(\Gamma, \varepsilon)$. 

It is not difficult to check that, both in the case $m=0$ (where we define $\bar{\mathcal  T} = \mathcal  T(\Gamma, \varepsilon)$ for sake of notational simplicity) and in the case $m>0$, the trisection $\bar{\mathcal  T}$  is completely identified by the attachments on $\Sigma$ of $H_{01}=H_0 \cap H_1=\bar H_0 \cap H_{1}$, $H_{02}=H_0 \cap H_2=\bar H_0 \cap H_{2}$ and $H_{12}=H_1 \cap H_2$ (which are 3-dimensional handlebodies by hypothesis).  
Now, as pointed out by Martini and Toriumi\footnote{See \cite[Paragraph 4.6]{Martini-Toriumi}, and in particular the caption of Figure 17.}, the quadrangulation of the central surface $\Sigma$ given by the union of the ``orange" squares in Figure \ref{fig:trisec_pieces} is dual to the cellular subdivision of the embedding surface $F_\e(\Gamma_{\hat 4})$, given by $\Gamma_{\hat 4}$ itself. Hence, $\Sigma$ and  $F_\e(\Gamma_{\hat 4})$ can be identified. 

Moreover, each $\{\e_0,\e_2\}$-colored (resp. $\{\e_1,\e_3\}$-colored) cycle of $\Gamma$, that is embedded in $F_\e(\Gamma_{\hat 4})$, is dual in $K(\Gamma_{\hat 4})$ to a red (resp. green) edge to which a prism constituting $H_{02}$ (resp. $H_{01}$) collapses (see Figure \ref{fig:trisec_pieces} and  \cite[Paragraph 4.6]{Martini-Toriumi}); further, the union of such edges, that is a spine of $H_{02}$ (resp. $H_{01}$), is isomorphic to $K_{\e_1 \e_3}$ (resp. $K_{\e_0 \e_2}$).
Hence, the statement regarding $\alpha$ and $\beta$ curves easily follows, by noting that a system of independent meridian curves for $H_{02}$ (resp. $H_{01}$) can be obtained from the set of  $\{\e_0,\e_2\}$-colored cycles (resp. $\{\e_1,\e_3\}$-colored cycles) by removing those corresponding to a maximal tree of $K_{\e_1 \e_3}$ (resp. $K_{\e_0 \e_2}$).

As regards $H_{12},$ note that it always collapses to the 2-dimensional complex $Q(\G, \e)$, consisting of a square for each $4$-colored edge of $\Gamma$, since in the cubical intersection of $H_{12}$ with any $4$-simplex the face belonging to the central surface is free:  see the proof of Proposition \ref{CS gem-induced trisections}, together with Figures \ref{fig:trisec_pieces} and \ref{figHsquare}. 

Moreover, by hypothesis, $Q(\G, \e)$ further collapses to a 1-dimensional complex via the collapsing sequence given by the ordering $(e_1, \dots, e_p)$ of $4$-colored edges satisfying condition (*), as described in the proof of Proposition \ref{CS gem-induced trisections}.  
Hence, a spine of $H_{12}$ consists of the set of edges of $Q(\G, \e)$ remaining at the end of the collapsing sequence, and a system of (independent) curves realizing the  attachment  of $H_{12}$ can  be obtained from this set of edges by shrinking to a point the ones belonging to a maximal tree.   
 The statement regarding the $\gamma$ curves follows, by projecting each element of the above system of curves (corresponding to a $\{c,4\}$-cycle of $\G$) on $F_\e(\Gamma_{\hat 4})$: in fact, while all $c$-colored edges ($c \in \Delta_3$) obviously embed into $F_\e(\Gamma_{\hat 4})$, the ordering (*) enables to project on $F_\e(\Gamma_{\hat 4})$ also the 4-colored edge $e_j$, for increasing values of $j\in\{1,\ldots,p\}$, via the bicolored cycle containing it whose edges different from $e_j$ already embed in $F_\e(\Gamma_{\hat 4})$.\qed
\smallskip

\begin{example}
{\em Figures \ref{fig:S^xS^2}, \ref{S2xS2_cycles_trisection-diagram} and \ref{fig:S^xS^2_trisection-diagram} show how to apply Proposition \ref{trisection_diagrams_gem-induced}: the gem of $\mathbb S^2 \times \mathbb S^2$ in Figure \ref{fig:S^xS^2} admits the depicted ordering of its $4$-colored edges, $(e_1, e_2, e_3, e_4, e_5, e_6, e_7)$,  which corresponds to a sequence of collapses of the squares of $Q(\G, \e)$ 
(where $\e=(0,1,2,3,4)$) from their free edges corresponding respectively to $\{c_i,4\}$-colored cycles, $\forall i = 1, \dots, 7$ (with $c_1=1$, $c_2=1$, $c_3=3$, $c_4=2$, $c_5=3$, $c_6=2$, $c_7=1$). The result is a 1-dimensional complex  consisting only of five edges: three corresponding to all $\{0,4\}$-colored cycles of $\Gamma$, one corresponding to the $\{2,4\}$-colored cycle containing $\{e_1, e_2, e_5, e_7\}$ and one corresponding to the $\{3,4\}$-colored cycle containing $\{e_2, e_4, e_6, e_7\}$. 

By shrinking to a point one of the edges corresponding to a $\{0,4\}$-colored cycle (for example the one containing $e_2$, $e_5$ and $e_4$), and both edges corresponding to $\{2,4\}$- and $\{3,4\}$-colored cycles, the system $\gamma$ of curves highlighted in red in Figure\ref{S2xS2_cycles_trisection-diagram}
is identified; it gives rise to a trisection diagram of $\mathbb S^2 \times \mathbb S^2$, together with the system $\alpha$ of curves highlighted in green (all $\{0,2\}$-colored cycles of $\Gamma$, but one arbitrarily chosen, since $K_{\e_1 \e_3}$ contains exactly two vertices) and the system $\beta$ of curves highlighted in blue (all $\{1,3\}$-colored cycles of $\Gamma$, but one arbitrarily chosen, since $K_{\e_0 \e_2}$ contains exactly two vertices).  
Figure  \ref{fig:S^xS^2_trisection-diagram} shows the associated trisection diagram of $\mathbb S^2 \times \mathbb S^2$.

Note that the three systems of curves on the (genus two) surface where $\G_{\hat 4}$ regularly embeds, in pairs give rise to genus two Heegaard diagrams of the boundaries of the 4-dimensional pieces of the trisection of $\mathbb S^2 \times \mathbb S^2$ (which are actually $3$-spheres).  }

\begin{figure} [ht]  
    \centering
    \includegraphics[width=0.8\linewidth]{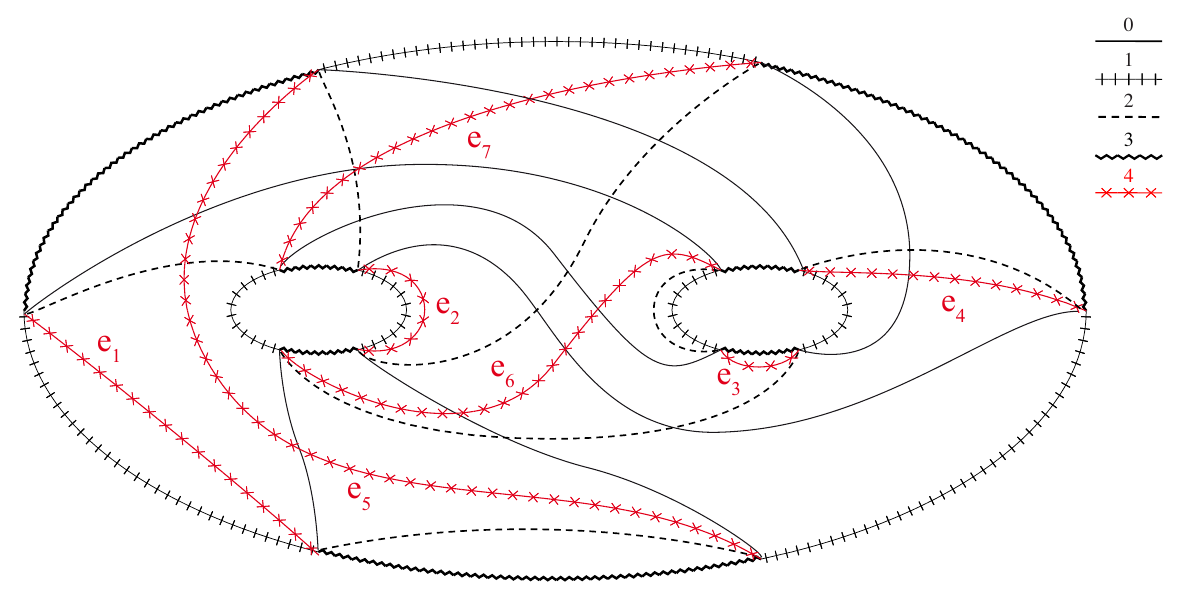}
   \caption{A gem of $\mathbb S^2 \times \mathbb S^2$, with an ordering of $4$-colored edges giving rise to a gem-induced trisection}
   \label{fig:S^xS^2}
\end{figure}

\begin{figure} [ht]  
    \centering
    \includegraphics[width=0.8\linewidth]{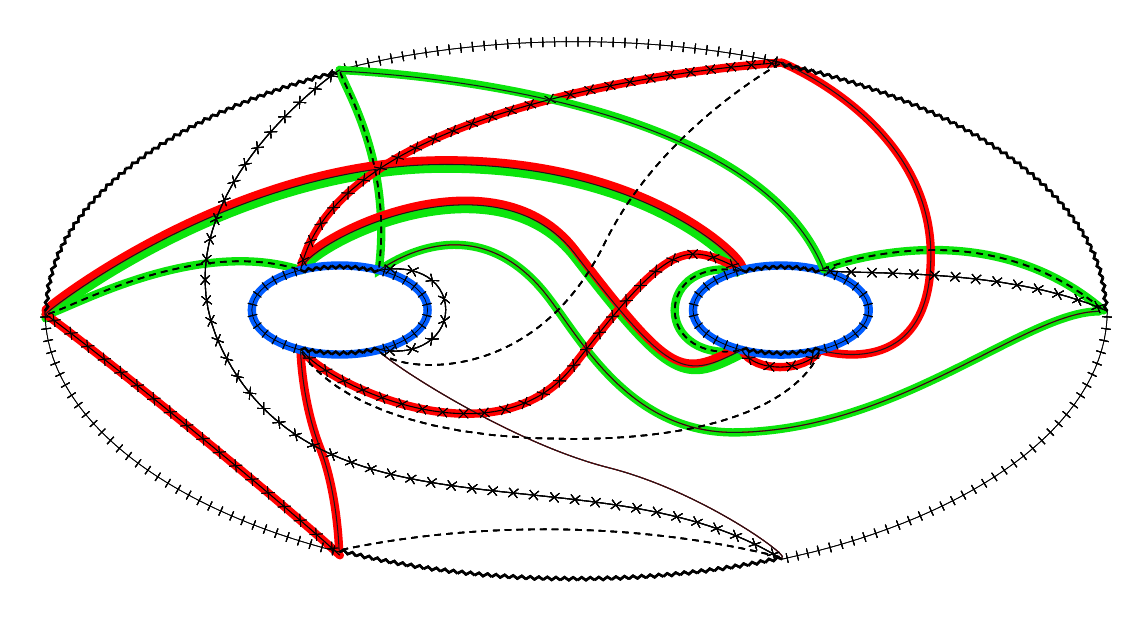}
   \caption{The bicolored cycles 
   giving rise to a trisection diagram}
   \label{S2xS2_cycles_trisection-diagram}
\end{figure}

\begin{figure} [!ht] 
    \centering
    \includegraphics[width=0.8\linewidth]{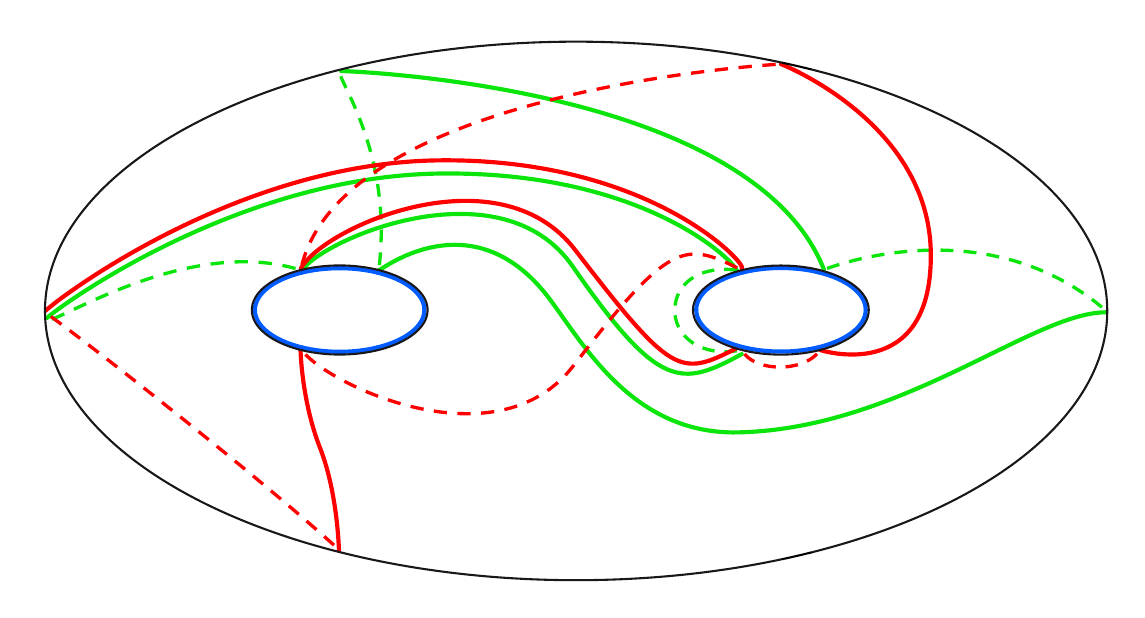}
    \caption{The associated trisection diagram  of $\mathbb S^2 \times \mathbb S^2$}
   \label{fig:S^xS^2_trisection-diagram}
\end{figure}

\end{example}

\bigskip
Let us now consider a (connected) Kirby diagram $(L,d)$, whose dotted components - if any - are in ``good position", which represents a closed 4-manifold $\bar M$. If $M$ is the compact $4$-manifold represented by the same diagram, $\bar M$ turns out to be obtained from $M$ by adding a suitable number of 3-handles. 
Since the gem $\Gamma(L,d)$ of $M$ satisfies condition (*) of Proposition \ref{CS gem-induced trisections} (see the proof of Theorem \ref{th_M4(L,c)}), the procedure described in Proposition \ref{trisection_diagrams_gem-induced} can always be applied in order to yield a trisection diagram for $\bar M$.  

\medskip 
Moreover, if no dotted component appears (i.e., if  $(L,d)$ is a framed link) the procedure described in Proposition \ref{trisection_diagrams_gem-induced}  can be applied to both the gem $\Gamma(L,d)$ of $M$ and the gem $\tilde \Omega(L,d)$ obtained from $\Gamma(L,d)$ via generalized dipole eliminations (see Section \ref{ss.from_Kirby_to_gems}). 

\smallskip 
The following statement proves that, if $(L,d)$ is a (connected) framed link representing a closed $4$-manifold $\bar M$, the three systems of curves of the associated trisection diagram are drawn on the closed orientable surface $\Sigma$  where $\Gamma_{\hat 4}(L,d)$ (resp. $\tilde \Omega_{\hat 4}(L,d)$) regularly embeds, whose genus is $s+1$, $s$ being the crossing number of $L$ (resp. is $m_\alpha$, where $m_\alpha$  is the number of $\alpha$-colored regions in a chess-board coloration of $L$, $\alpha$ being the color of the unbounded region), according to Proposition \ref{order-regular genus Gamma}.

Moreover, each system of curves is simply obtained from the cycles of a suitable pair of colors of $\Gamma(L,d)$ (resp. $\tilde\Omega(L,d)$): 

\begin{proposition} \label{trisection_diagrams_framed-links}
With the above notations, let $\mathcal T$ be the trisection of $\bar M$ induced by $\Gamma(L,d)$ (resp. $\tilde \Omega(L,d)$). Then, a trisection diagram for $\mathcal  T$ is given by the following three system of curves $\{ \alpha, \beta, \gamma\}$ on $\Sigma$:

\begin{itemize}
\item[-]   the $\alpha$ curves are all $\{1,2\}$-colored cycles of $\Gamma(L,d)$ (resp. $\tilde\Omega(L,d)$), but one arbitrarily chosen;   
\item[-] the $\beta$ curves are all $\{0,3\}$-colored cycles of $\Gamma(L,d)$ (resp. $\tilde\Omega(L,d)$),  but $2l-1$ corresponding to a maximal tree of the 1-dimensional subcomplex $K_{12}$ generated by the vertices labelled $1$ and $2$; 

\item[-] the $\gamma$ curves are all $\{2,4\}$-colored cycles of $\Gamma(L,d)$ (resp. $\tilde\Omega(L,d)$), but one arbitrarily chosen. 
\end{itemize}
\end{proposition} 

\dimo
Proposition \ref{order-regular genus Gamma} ensures that, if the permutation $\varepsilon = (1,0,2,3,4)$ of $\Delta_4$ is chosen, then $\rho_{\varepsilon}(\G_{\hat 4}(L,d))=\rho_{\varepsilon}(\Lambda(L,d)) = s + 1$ (resp. $\rho_{\varepsilon}(\tilde\Omega_{\hat 4}(L,d))=\rho_{\varepsilon}(\Omega(L,d)) = \  m_\alpha$). 
Moreover, by construction, all 4-colored edges of $\Gamma(L,d)$ (resp.  $\tilde \Omega(L,d)$) double $1$-colored edges of $\Gamma_{\hat 4}(L,d)= \Lambda(L,c)$ (resp. of  $\tilde \Omega_{\hat 4}(L,d) =\Omega(L,d)$), but the $3l$ ones added between the vertices $P_{2r}^{(i)}$ and $P_{2r+1}^{(i)}$, for $r=0,1,2$, of the quadricolor corresponding to the $i$-th component of $L$, $\forall i \in \{1, 2, \dots, l\}$.

Let us now consider the following ordering of all $4$-colored edges $\Gamma(L,d)$ (resp.  $\tilde \Omega(L,d)$): 
\begin{itemize}
    \item[-] first, in any order, all $4$-colored edges belonging to $\{1,4\}$-colored cycles of length two, together with the $l$ $4$-colored edges between the vertices $P_0^{(i)}$ and $P_{1}^{(i)}$, $\forall i \in \{1, 2, \dots, l\}$ (which belong to $\{0,4\}$-colored cycles of length two);
    \item[-] then, the $l$ $4$-colored edges between the vertices $P_4^{(i)}$ and $P_{5}^{(i)}$, $\forall i \in \{1, 2, \dots, l\};$ 
        \item[-] finally, the $l$ $4$-colored edges between the vertices $P_2^{(i)}$ and $P_{3}^{(i)}$, $\forall i \in \{1, 2, \dots, l\}.$ 
\end{itemize}
It is not difficult to check that the above ordering - which satisfies condition (*) of Proposition  \ref{CS gem-induced trisections} - corresponds to a collapsing sequence of all squares of the 2-dimensional complex $Q(\G, \e)$, giving rise to a 1-dimensional spine  
$Q^\prime$ of  $H_{12}$ which consists of the edges corresponding to: 
\begin{itemize}
    \item[-] all $\{2,4\}$-colored cycles of $\Gamma(L,d)$ (resp.  $\tilde \Omega(L,d)$) 
    \item[-]  all $\{0,4\}$-colored cycles of $\Gamma(L,d)$ (resp.  $\tilde \Omega(L,d)$), but those containing  $P_0^{(i)}$ and $P_{1}^{(i)}$, $\forall i \in \{1, 2, \dots, l\}$
    \item[-]  all $\{3,4\}$-colored cycles of $\Gamma(L,d)$ (resp.  $\tilde \Omega(L,d)$), but those containing  $P_4^{(i)}$ and $P_{5}^{(i)}$, $\forall i \in \{1, 2, \dots, l\}.$  
\end{itemize}

Moreover, the combinatorial structure of $\Gamma(L,d)$ and  $\tilde \Omega(L,d)$ ensures that all edges of $Q^\prime$  corresponding to $\{0,4\}$-colored (resp. $\{3,4\}$-colored) cycles have a free end-point, since the corresponding cycles belong to disjoint $\{0,1,4\}$-residues (resp. $\{1,3,4\}$-residues). Hence, $Q^\prime$  further collapses to a 1-dimensional subcomplex, whose edges correspond to the $\{2,4\}$-colored cycles of $\Gamma(L,d)$ (resp.  $\tilde \Omega(L,d)$).       

On the other hand, all these edges have the same end-points, since  $\Gamma(L,d)$ (resp.  $\tilde \Omega(L,d)$) contains exactly one $\{0,2,4\}$-residue and exactly one $\{2,3,4\}$-residue (see the proof of \cite[Theorem 2]{Casali JKTR2000}, ensuring that $\Lambda_{\hat j}(L,d)$, and hence $\Omega_{\hat j}(L,d)$, is connected, for $j \in \{0,3\}$). 
The statement concerning the $\gamma$ curves now directly follows.

Finally, the statement concerning the $\alpha$ (resp. $\beta$) curves is a consequence of the analogue general statement of Proposition \ref{trisection_diagrams_gem-induced}, together with the fact that both $\Gamma_{\hat 4}(L,d)=\Lambda(L,d)$ and $\tilde \Omega_{\hat 4}(L,d)=\Omega(L,d)$ have exactly one $\hat 0$-residue and one $\hat 3$-residue  (resp. exactly $l$ $\hat 1$-residue and $l$ $\hat 2$-residue), according to the proof of \cite[Theorem 2]{Casali JKTR2000}. \qed

\medskip

Finally, we point out that our particular extension to the boundary case of the notion of trisection, given by gem-induced trisections, suggests also a possible extension of the notion of trisection diagram to simply-connected compact $4$-manifolds with connected boundary. 

\begin{definition}
{\em A {\it G-trisection diagram}  of genus $g$ is a 4-tuple $(\Sigma; \alpha, \beta, \gamma)$, where $\Sigma$ is a genus $g$ orientable surface and $\alpha$, $\beta$, $\gamma$ are complete systems of curves on $\Sigma$, such that: 
\begin{itemize}
    \item[-] 
  $(\Sigma; \alpha, \gamma)$ and $(\Sigma; \beta, \gamma)$ are (genus $g$) Heegaard diagrams of $\#_{k_1}(\mathbb S^1 \times \mathbb S^2)$ and $\#_{k_2}(\mathbb S^1 \times \mathbb S^2)$ respectively, with $0\le k_i \le g$, for $i=1,2$;
  \item[-] 
 $(\Sigma; \alpha, \beta)$ is a (genus $g$) Heegaard diagram of a closed connected 3-manifold.  
\end{itemize} 
 }
\end{definition}

 \begin{proposition}   \label{trisection_diagrams_simply-connected_bounded}
     Let $M$ be a simply-connected compact $4$-manifold with connected boundary, admitting a genus $g$ gem-induced trisection. 
     Then, $M$ is uniquely identified by a G-trisection diagram $(\Sigma; \alpha, \beta, \gamma)$ of genus $g$, 
     so that $(\Sigma; \alpha, \beta)$ is a Heegaard diagram of $\partial M$.    
 \end{proposition} 

 \dimo
Let $\mathcal  T = (H_{0},H_{1},H_{2})$ be a triple of submanifolds of $M$ constituting a gem-induced trisection. By definition,  $H_{1}$ and $H_{2}$ are $4$-dimensional handlebodies, while $H_{0}$ is a collar of $\partial M$; moreover, all pairwise intersections are 3-dimensional handlebodies, while $\Sigma = H_{0} \cap H_{1} \cap H_{2}$ is a closed orientable genus $g$ surface. 
If $\alpha, \beta, \gamma$ are attaching curves on $\Sigma$ of $H_{01}=H_0 \cap H_1$, $H_{02} = H_0 \cap H_2$ and $H_{12}=H_1 \cap H_2$ respectively,  
then $(\Sigma; \alpha, \beta, \gamma)$ turns out to be a G-trisection diagram (see \cite[Remark 15]{Casali-Cristofori gem-induced} for details). 

On the other hand, given such a G-trisection diagram, let $W$ be the compact 4-manifold obtained by attaching a collar of two genus $g$ 3-dimensional handlebodies  (corresponding to $H_{01}$ and $H_{02}$)  to $\Sigma \times \mathbb D^2$  according to the curves $\alpha$ and $\beta$ respectively, in full analogy with the procedure described after Definition \ref{def. trisection-diagram}.  Since $M$ is obtained from $W$ by adding 3-handles (corresponding to $H_1$ and $H_2$), and $M$ is simply-connected with connected boundary by hypothesis, Theorem 1 in \cite{[Trace 1982]}  
ensures that $W$ (or, equivalently, the G-trisection diagram) uniquely determines $M$.\qed

\medskip

As a consequence, the following statement easily follows from Proposition \ref{trisection_diagrams_gem-induced}, by adapting the same proving arguments to the (general) boundary case: 

    \begin{proposition} \label{trisection_diagrams_gem-induced_boundary} \ 
Let $\Gamma\in G_s^{(4)}$ be a gem of a simply-connected compact $4$-manifold $M$ with connected boundary $\partial M$, satisfying condition (*) of Proposition \ref{CS gem-induced trisections}. 

Then, for any cyclic permutation $\varepsilon$ of $\Delta_4$, a G-trisection diagram for $M$ is given by the following three system of curves $\{ \alpha, \beta, \gamma\}$ on $F_\e(\Gamma_{\hat 4})$:
\begin{itemize}
\item[-] the $\alpha$ curves (resp. $\beta$ curves) consist of all $\{\e_0,\e_2\}$-colored cycles (resp. $\{\e_1,\e_3\}$-colored cycles) of $\Gamma$, but those corresponding to a maximal tree of the 1-dimensional subcomplex  
$K_{\e_1 \e_3}$ (resp. $K_{\e_0 \e_2}$),  generated by the vertices labelled $\e_1$ and $\e_3$ (resp. $\e_0$ and $\e_2$) of $K(\Gamma).$    
\item[-] the $\gamma$ curves are a subset of  $\{c,4\}$-colored cycles of $\Gamma$, with $c \in \Delta_3,$ which can be combinatorially identified by following - via condition (*) - the collapsing sequence for the 3-dimensional handlebody $H_{12}$ (according to the notations of Section \ref{sss Gem-induced trisections}).   
\end{itemize} 
\end{proposition} 

\bigskip 

Also Proposition \ref{trisection_diagrams_framed-links} can be extended to the case of (connected) framed links representing $4$-manifolds with boundary:

 \begin{proposition}  \label{trisection_diagrams_gem-induced_framed-links}
Let $(L,d)$ be a (connected) framed link representing a (simply-connected) $4$-manifold $M$ with non-empty boundary, and let $\mathcal T$ be the gem-induced trisection of $M$ induced by $\Gamma(L,d)$ (resp. $\tilde \Omega(L,d)$). 
Then, a G-trisection diagram for $\mathcal  T$ is given by the following three system of curves $\{ \alpha, \beta, \gamma\}$ on $\Sigma$:
\begin{itemize}
    \item[-]  the $\alpha$ curves are all $\{1,2\}$-colored cycles of $\Gamma(L,d)$ (resp. $\tilde\Omega(L,d)$),
but one arbitrarily chosen;   

\item[-]  the $\beta$ curves are all $\{0,3\}$-colored cycles of $\Gamma(L,d)$ (resp. $\tilde\Omega(L,d)$), 
but $2l-1$ corresponding to a maximal tree of the 1-dimensional subcomplex $K_{12}$    generated by the vertices labelled $1$ and $2$;

\item[-]  the $\gamma$ curves are all $\{2,4\}$-colored cycles of $\Gamma(L,d)$ (resp. $\tilde\Omega(L,d)$), but one arbitrarily chosen.    
\end{itemize} 
\end{proposition} 

\smallskip

\begin{remark} \label{meridians_as_gamma_curves} {\em It is not difficult to prove that, both in Proposition \ref{trisection_diagrams_framed-links} and in Proposition \ref{trisection_diagrams_gem-induced_framed-links}, a slight modification of the ordering of $4$-colored edges (by considering the $l$ $4$-colored edges between $P_2^{(i)}$ and $P_{3}^{(i)}$, $\forall i \in \{1, 2, \dots, l\}$, before the $l$ ones between $P_4^{(i)}$ and $P_{5}^{(i)}$, $\forall i \in \{1, 2, \dots, l\}$) 
yields another G-trisection  
diagram associated to $\Gamma(L,d)$ (resp. $\tilde \Omega(L,d)$), with the same $\alpha$-curves and $\beta$-curves, while the $\gamma$ curves consist of: 
\begin{itemize}
 \item[$\bullet$] $l$ curves corresponding to the $\{3,4\}-$cycles of $\Gamma(L,d)$ (resp. $\tilde\Omega(L,d)$) containing the vertices $P_4^{(i)}$ and $P_{5}^{(i)}$, $\forall i \in \{1, 2, \dots, l\}$; 
 \item[$\bullet$] $s+1-l$ (resp. $m_\alpha -l$) curves corresponding to all $\{2,4\}-$cycles of $\Gamma(L,d)$ (resp. $\tilde\Omega(L,d)$) different from the length two ones with vertices $P_2^{(i)}$ and $P_{3}^{(i)}$, $\forall i \in \{1, 2, \dots, l\}$, but one arbitrarily chosen.
\end{itemize} 
}
\end{remark}

\medskip

\begin{example}{\em
Figures \ref{fig:M4(trefoil,1)}, \ref{fig:M4(trefoil,1)_cycles_trisection-diagram} and \ref{fig:M4(trefoil,1)_trisection-diagram} show an application of Proposition \ref{trisection_diagrams_gem-induced_framed-links}, also making use of Remark \ref{meridians_as_gamma_curves}, to the framed link consisting of a trefoil knot $K_T$ with framing 1: the gem $\tilde\Omega(K_T,1)$ of Figure \ref{fig:M4(trefoil,1)} represents a (simply-connected) compact $4$-manifold $M$ (whose boundary is the Poincar\'e homology sphere), has regular genus $m_\alpha +1 = 3$ and gives rise to a gem-induced trisection of genus $m_\alpha=2$. The bicolored cycles giving rise to the associated G-trisection diagram are highlighted in Figure \ref{fig:M4(trefoil,1)_cycles_trisection-diagram}, while the associated G-trisection diagram is depicted in Figure \ref{fig:M4(trefoil,1)_trisection-diagram}. 
Note that by considering blue and red  
(resp. red 
and green) curves, a genus two Heegaard diagram of the boundary of one of the 4-dimensional handlebodies 
 of the trisection  of $M$ is obtained, while blue 
and green curves give rise to a genus two Heegaard diagram of the Poincar\'e homology sphere. }
\end{example}

\begin{figure} [!ht]  
    \centering
    \includegraphics[width=0.6\linewidth]{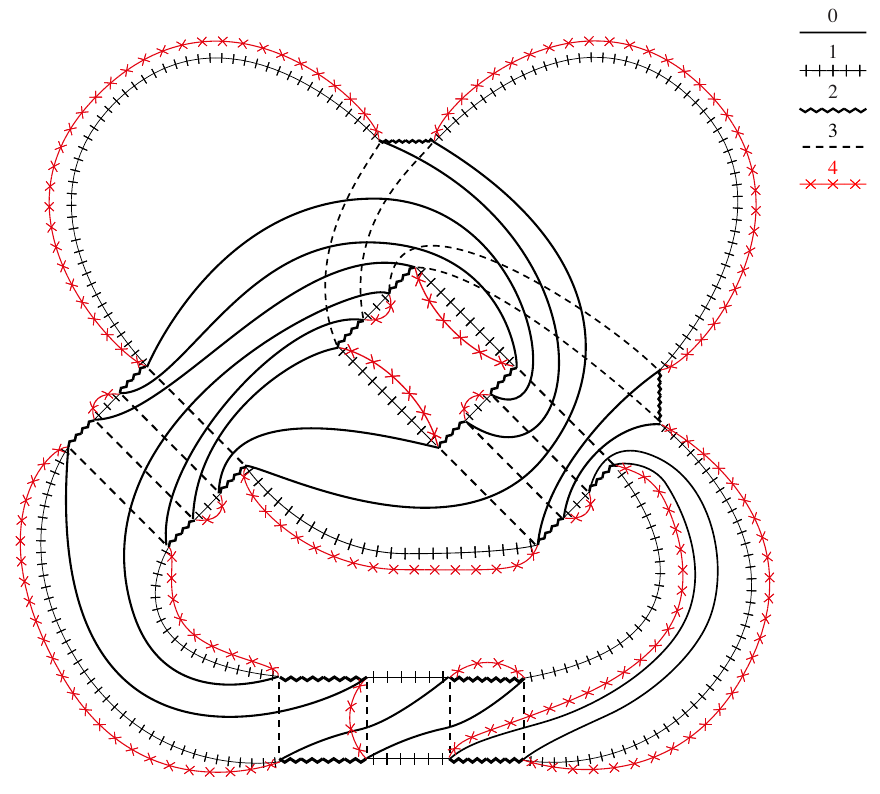}
    \caption{The gem $\tilde\Omega(K_T,1)$, $K_T$ being the trefoil knot}
    \label{fig:M4(trefoil,1)} 
\end{figure}

\begin{figure} [!ht]  
    \centering
    \includegraphics[width=0.6\linewidth]{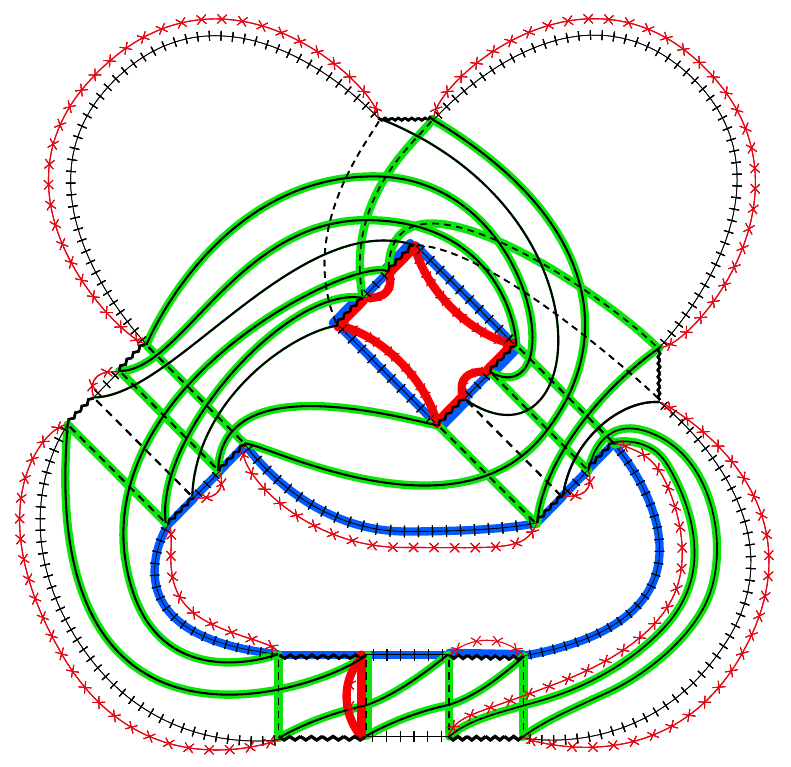}
    \caption{The cycles of $\tilde\Omega(K_T,1)$, giving rise to a G-trisection diagram}
    \label{fig:M4(trefoil,1)_cycles_trisection-diagram} 
\end{figure}

 \begin{figure} [!t] 
     \centering
     \includegraphics[width=0.6\linewidth]{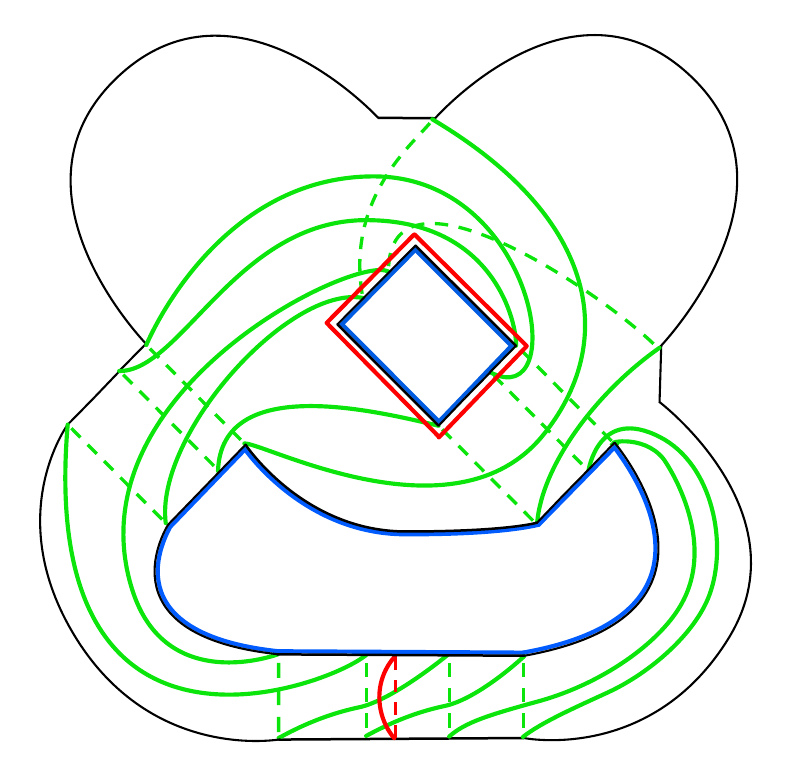}
     \caption{The G-trisection diagram arising from $\tilde\Omega(K_T,1)$} 
     \label{fig:M4(trefoil,1)_trisection-diagram}
 \end{figure}

\bigskip

\begin{figure} [t] 
    \centering
    \includegraphics[width=0.6\linewidth]{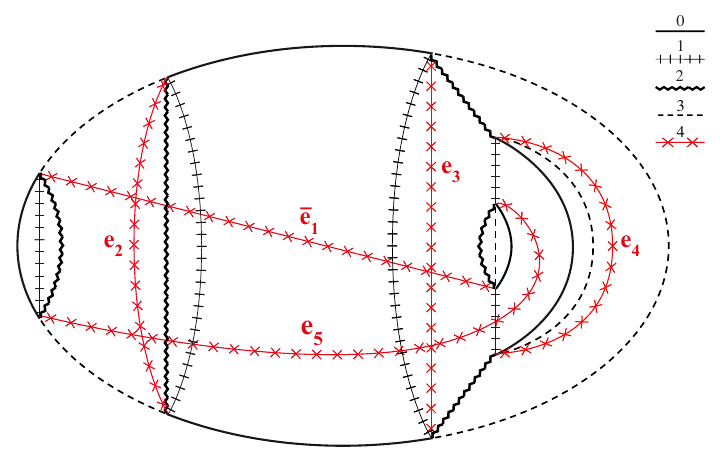}
 \caption{the standard order ten gem of $\mathbb S^1\times\mathbb S^3$ with the stabilizing edge $\bar e_1$}
    \label{fig:S3xS1_stabilization} 
\end{figure}

\begin{figure} [!ht] 
    \centering
    \includegraphics[width=0.6\linewidth]{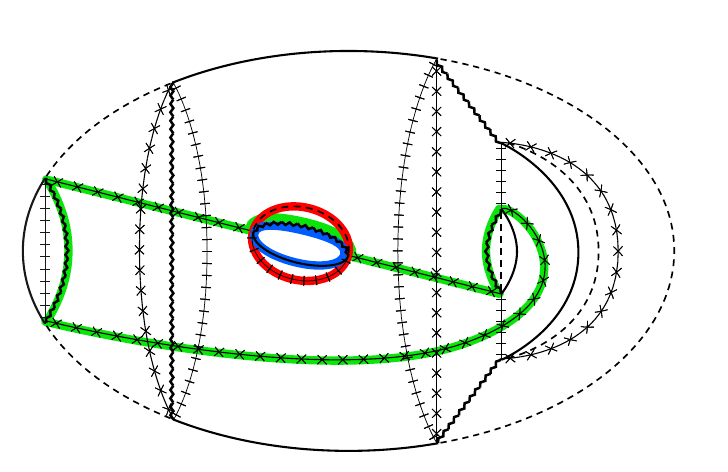}    
 \caption{the curves of the trisection diagram of  $\mathbb S^1\times\mathbb S^3$ induced by the gem}
    \label{fig:S3xS1_curves} 
\end{figure}

\begin{figure} [!ht] 
    \centering
     \includegraphics[width=0.3\linewidth]{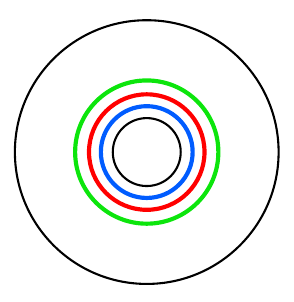}
\caption{the associated trisection diagram of $\mathbb S^1\times\mathbb S^3$}
    \label{fig:S3xS1_standard_diagram} 
\end{figure}

 \begin{remark}{\em In \cite{Casali-Cristofori trisection ter}, by using  stabilizations along 4-colored edges and the ordering  involved in the proof of Theorem \ref{number_stabilizations}, both Propositions \ref{trisection_diagrams_gem-induced} and \ref{trisection_diagrams_simply-connected_bounded} are extended to gems that do not directly induce trisections. As a consequence, trisection diagrams can be obtained algorithmically from \underline{any}  gem of a closed 4-manifold.  
Detailed statements and proofs can be found in \cite{Casali-Cristofori trisection ter}, while Figures \ref{fig:S3xS1_stabilization}, \ref{fig:S3xS1_curves} and \ref{fig:S3xS1_standard_diagram}  show an example of the procedure, applied to the standard order ten crystallization $\G$ of $\mathbb S^1 \times \mathbb S^3$ depicted in Figure \ref{fig:S3xS1_stabilization}. If a stabilization along the 4-colored edge $\bar e_1$ is performed, all squares of $Q(\G, \e)$ ($\e=(0,1,2,3,4)$) subsequently collapse and the obtained 1-dimensional complex, after shrinking to a point a maximal tree, corresponds to the $\{2,4\}$-colored cycle highlighted in red in Figure \ref{fig:S3xS1_curves}. Hence, this cycle identifies the (unique) curve of the system $\gamma$ in the associated trisection diagram of genus one, whose surface is obtained, via the stabilization, from the 2-sphere where $\G_{\hat 4}$ regularly embeds. On the other hand, the systems $\alpha$ and $\beta$ consist of the blue and green trivial curves respectively, which are meridian curves of the 1-handles involved in the stabilization. It is easy to check that the obtained trisection diagram exactly coincides with the standard one of $\mathbb S^1 \times \mathbb S^3$ (see Figure \ref{fig:S3xS1_standard_diagram}).
 } \end{remark}

\bigskip

\section{Open problems and trends}

The constructions and results presented in the previous Sections show the strong relationship which exists among the three involved representations for PL 4-manifolds, namely Kirby diagrams, trisections and gems, pointing out how the latter can be useful in connecting the first two, also in order to compute or estimate invariants.      

On the other hand, the analysis carried out also suggests several possible developments, which could help to complete the picture of the existing relationships among the above theories. Our hope is that a unifying vision can offer deep and unexpected insights, allowing to improve the understanding of the topic.  

\medskip

As regards Section \ref{ss.from_Kirby_to_gems}, an interesting development would be to obtain an algorithmic procedure to pass from Kirby diagrams to gems
in the generalized case of diagrams representing also non-orientable 1-handles (according to \cite{Caesar De Sa}). 
If the resulting gems admitted gem-induced trisections (as happens in the orientable setting: see Theorem \ref{trisection_vs_regular-genus(closed case)}), in virtue of the non-orientable version of Laudenbach-Poenaru's result due to \cite{Miller-Naylor}, an estimation of the trisection genus for \underline{all} closed 4-manifolds would be obtained, via gems. 

\smallskip 

Another possible improvement would be to study how to obtain a gem of a closed 4-manifold, starting from a gem of the compact one consisting only of handles up to index two, even if no $\rho_3$- or $\rho_4$-pair appears in the latter (as assumed in Proposition \ref{chiudere}). 
The problem is obviously a difficult one, since it involves the combinatorial recognition of $\#_m(\mathbb S^1 \times \mathbb S^2)$ ($m >0$) in 3-dimensional topology, but in specific situations it could possibly be faced, for example by detecting sequences of combinatorial moves yielding the desired $\rho_3$-pairs. 

A particular case of this problem arises, for instance, in the construction of a gem for the double of a compact $4$-manifold. It is known that a Kirby diagram of the double can be obtained,  from a Kirby diagram of the manifold itself, by adding a 0-framed meridian to each framed component (see \cite{GS}).  
Therefore,  it may be worth investigating whether the particular configuration of the arising 5-colored graph allows to be modified in a standard way so as to give rise to $\rho_3$-pairs satisfying Proposition  \ref{chiudere}.

\medskip

With respect to Section \ref{ss.from_gems_to_trisections}, a fundamental problem involves the possible characterization of classes of gems so that the arising trisections (either directly, as gem-induced trisections, or indirectly, via the approach described in Section \ref{ss: First indirect approach}) realize the trisection genus of the represented $4$-manifolds. 

In \cite[Section 6]{Casali-Cristofori gem-induced}, a detailed analysis has already been performed in the case of $4$-manifolds (with empty or connected boundary) admitting a gem-induced trisection $\mathcal T(\Gamma, \e)$ so that $\Gamma$ is {\it weak semi-simple}\footnote{A $5$-colored graph $\G\in G_s^{(4)}$,  representing a compact $4$-manifold $M$ with empty or connected boundary, is said to be {\em weak  semi-simple} with respect to a cyclic permutation $\varepsilon$ of $\Delta_4$ if \, $g_{\varepsilon_i,\varepsilon_{i+2},\varepsilon_{i+4}} = 1 + m \ \ \forall \ i \in \{0,2,4\}$ \, and \, $g_{\varepsilon_i,\varepsilon_{i+2},\varepsilon_{i+4}} = 1 + m^{\prime} \ \ \forall \ i \in \{1,3\}$ \,  where  \, $rk(\pi_1(M))= m \ge 0$ and  $rk(\pi_1(\widehat M))= m^{\prime} \ge 0$.} with respect to the cyclic permutation $\e$ of the color set. 
As a consequence, a large class of closed simply-connected 4-manifolds - comprehending all standard ones - has been identified for which the trisection genus equals the second Betti number and also coincides with half the regular genus, and an even larger class of such manifolds, for which the trisection genus equals the second Betti number\footnote{Note that for both classes the involved minimal trisections are {\it efficient}, according to \cite{Lambert-Cole-Meier}.}:   
\medskip

\begin{proposition} \label{minimality closed weak simple}  \  {\rm (\cite{Casali-Cristofori gem-induced})} \
$$g_T(M)\ =\ \beta_2(M)$$
for each closed  (simply-connected) $4$-manifold $M$  which admits a $5$-colored graph $\G\in G_s^{(4)}$ representing $M$ with  
$g_{\varepsilon_0,\varepsilon_{2},\varepsilon_{4}} = 1 $ and 
$g_{\varepsilon_1,\varepsilon_{3},\varepsilon_{4}} = 1 $ 
($\varepsilon$ being a cyclic permutation of $\Delta_4$), so that $\mathcal T(\G,\varepsilon)$ is a gem-induced trisection. 
\par \noindent 
In particular: 
$$g_T(M)\ =\ \beta_2(M) =\frac 1 2 
 \, \mathcal G(M)$$      
for each closed  (simply-connected) $4$-manifold $M$ admitting a  crystallization $\G$ with $g_{\varepsilon_i,\varepsilon_{i+2},\varepsilon_{i+4}} = 1 \ \ \forall i \in \Delta_4$ ($\varepsilon$ being a cyclic permutation of $\Delta_4$)\footnote{This condition is usually referred to by saying that  $\G$ is {\it weak simple} with respect to $\varepsilon$.} and such that $\mathcal T(\G,\varepsilon)$ is a gem-induced trisection.
\end{proposition}

It would be very interesting to characterize topologically $4$-manifolds admitting gems which satisfy the hypotheses of the above statement, or at least to find classes satisfying them.

On the other hand, it would also be interesting to investigate conditions on the gems of a compact $4$-manifold, so that in Theorem \ref{trisection_from_gem-induced} the equality holds.

\bigskip

Finally, we point out that some of the procedures described in Section \ref{ss.from_gems_to_trisections} could be performed by suitable automatic implementations. 
In particular: 
\begin{itemize}
    \item Starting from a given gem (for example, by using   its alphanumerical {\it code}: see \cite{Casali-Gagliardi code}), it is possible to determine algorithmically  
    the best choice of edges to use for stabilizations, so as to minimize the genus of the arising trisection, and hence get a better estimation of the trisection genus of the associated closed $4$-manifold.  
    \item Starting from a given gem, it is obviously possible to check in an automatic way whether condition (*) of Proposition \ref{CS gem-induced trisections} is satisfied; if so, it should not be difficult to detect algorithmically the bicolored cycles yielding a trisection diagram of the associated gem-induced trisection.
    \end{itemize}

\bigskip

\noindent {\bf Acknowledgements:} This work was supported by GNSAGA of INDAM and by the University of Modena and Reggio Emilia, project:  {\it ``Discrete Methods in Combinatorial Geometry and Geometric Topology''.}


\bigskip 


\begin{thebibliography}{00}


\bibitem{Bandieri-Gagliardi} 
P.~Bandieri, C.~Gagliardi: {\it Rigid gems in dimension $n$}, Bol. Soc. Mat. Mex. \textbf{18}(3), 55-67 (2012).   

\bibitem{Bell-et-al} 
M. Bell, J. Hass, J.~H. Rubinstein, S. Tillmann: {\it Computing trisections of 4-manifolds}, Proc. Nat. Acad. Sci. USA {\bf 115}(43), 10901-10907  (2018).

\bibitem{Burke} R.~A. Burke: {\it Practical Software for Triangulating and Simplifying 4-Manifolds}, Journal of Computational Geometry  {\bf  16}(2), 109-144 (2025).  

\bibitem{Burke-Burton-Spreer} R.~A. Burke – B. Burton – J. Spreer: {\it Small Triangulations of 4-Manifolds and the 4-Manifold Census}, Discrete Comput. Geom. (published online February 2026). \ \ https://doi.org/10.1007/s00454-026-00818-w


\bibitem{regina} B.A. Burton, R. Budney, W. Pettersson, et al.: {\it Regina: Software for low-dimensional topology}, http://regina-normal.github.io/, 1999–2025.    

\bibitem{Casali JKTR2000} 
M.~R.~Casali,  {\it From framed links to crystallizations of bounded 4-manifolds}, J. Knot Theory Ramifications {\bf 9} (4),  443-458  (2000).

\bibitem{Casali-Cristofori ElecJComb 2015} 
M.~R.~Casali, P.~Cristofori: {\it Cataloguing PL 4-manifolds by gem-complexity}, Electron. J. Combin. {\bf 22}(4), \# P4.25  (2015).

\bibitem{Casali-Cristofori Kirby-diagrams} 
M.~R. Casali, P. Cristofori: {\it Kirby diagrams and $5$-colored graphs representing compact $4$-manifolds}, Rev. Mat. Complut.  {\bf 36}, 899-931 (2023).   

\bibitem{generalized-genus}    
M.~R. Casali, P. Cristofori: {\it Classifying compact 4-manifolds via generalized regular genus and G-degree},  Ann. Inst. Henri Poincar\'e Comb. Phys. Interact.  {\bf 10}(1), 121-158 (2023). 

\bibitem{Casali-Cristofori gem-induced} 
M.~R. Casali, P. Cristofori: {\it Gem-induced trisections of compact PL 4-manifolds},  Forum Math. {\bf 36}(1), 87-109 (2024).  

\bibitem{Casali-Cristofori trisection bis} 
M.~R. Casali, P. Cristofori: {\it Trisections of PL 4-manifolds arising from colored triangulations}, Mediterr. J. Math. {\bf 22}, 27 (2025).  


\bibitem{Casali-Cristofori trisection ter} 
M.~R. Casali, P. Cristofori: {\it Estimating trisection genus via gem theory}, preprint (2025). \ \ DOI: 10.48550/arXiv.2504.04434


\bibitem{Casali-Cristofori-Gagliardi RIMUT 2020}
M.~R. Casali, P. Cristofori, C. Gagliardi:  {\it Crystallizations of compact 4-manifolds minimizing combinatorially defined PL-invariants}, Rend. Istit. Mat. Univ. Trieste {\bf 52}, 431-458 (2020).  

\bibitem{Casali-Cristofori-Grasselli}
M.~R. Casali, P. Cristofori, L. Grasselli: {\it G-degree for singular manifolds}, RACSAM  {\bf 112}(3), 693-704 (2018).  

\bibitem{Casali-Cristofori-Dartois-Grasselli}
M.~R. Casali, P. Cristofori,  S.~Dartois, L.~ Grasselli: {\it Topology in colored tensor models  via crystallization theory},  J. Geom. Phys. {\bf 129}, 142-167  (2018). 


\bibitem{Casali-Cristofori-Gagliardi Complutense 2015}
M.~R. Casali, P. Cristofori, C. Gagliardi: {\it Classifying PL 4-manifolds via crystallizations: results and open problems}, in:  {\em ``Mathematical Tribute to Professor Jos\'e Mar\'ia Montesinos Amilibia''}, Universidad Complutense Madrid (2016).  \ \ https://www.mat.ucm.es/~josefer/otros/homenaje-montesinos.pdf  

\bibitem{Casali-Gagliardi code}
M.~R. Casali, C. Gagliardi:  {\it A code for m-bipartite edge-coloured graphs}, Rend. Istit. Mat. Univ. Trieste {\bf 32},   55-76 (2001). 


 \bibitem{Castro-Gay-Pinzon} 
N.~A. Castro,  D.~T. Gay,  J. Pinzon-Caicedo: {\it Trisections of 4-manifolds with boundary},  Proc. Nat. Acad. Sci. USA  {\bf 115}(43), 10861-10868 (2018).  

\bibitem{Caesar De Sa} E.C. de S\'a: {\it A link calculus for 4-manifolds}. In: Fenn, R. (eds), \textit{Topology of Low-Dimensional Manifolds}, Lecture Notes in Mathematics {\bf 722}, Springer  (1979). 

\bibitem{Ferri-Gagliardi-Grasselli} 
M. Ferri, C. Gagliardi, L.  Grasselli: {\it A graph-theoretical representation of PL-manifolds. A survey on crystallizations}, Aequationes Math. {\bf 31}, 121-141 (1986).

\bibitem{Gagliardi 1987} C. Gagliardi, {\it On a class of 3-dimensional polyhedra}, Ann. Univ. Ferrara, Sez. VII, Sc. Mat., {\bf 33}, 51-88 (1987).

\bibitem{Gagliardi-Volzone}
C. Gagliardi, G. Volzone:  {\it Handles in Graphs and Sphere Bundles over $\mathbb S^1$},  Europ. J. Combinatorics  {\bf  8}, 151-158 (1987). 

\bibitem{Gay-Kirby}
D. Gay, R. Kirby:  {\it Trisecting 4-manifolds},  Geom. Topol.  {\bf 20}, 3097-3132 (2016).

\bibitem{GS} 
R.E. Gompf, A.I. Stipsicz:  {\it 4-manifolds and Kirby calculus}. American Mathematical Society, {\bf 20}, (1999).  

\bibitem{Lambert-Cole-Meier} 
P. Lambert-Cole, J. Meier: {\it Bridge trisections in rational surfaces},  Journal of Topology and Analysis {\bf 14}(3), 655-708 (2022). 

\bibitem{Laudenbach-Poenaru}
F. Laudenbach, V. Poenaru: {\it A note on 4-dimensional handlebodies}, Bull. Soc. Math. France  {\bf 100}, 337-344  (1972). 

\bibitem{Lins-book} S.~Lins: {\it Gems, computers and attractors for 3-manifolds}. Knots and Everything, no. 5. World Scientific, River Edge, NJ, 1995.

\bibitem{[M]} 
R. Mandelbaum:  {\it Four-dimensional topology: an introduction},  Bull. Amer. Math. Soc. {\bf 2}, 1-159 (1980).

\bibitem{Martini-Toriumi} R. Martini, R. Toriumi: {\it Trisections in colored tensor models},  Ann. Inst. Henri Poincar\'e Comb. Phys. Interact. {\bf 11}(3), 453–502 (2024).  

\bibitem{Meier} 
J. Meier: {\it Trisections and spun 4-manifolds},   Math. Res. Lett. {\bf 25}(5), 1497-1524 (2018).   

\bibitem{Meier-Schirmer-Zupan}
J. Meier, T. Schirmer, A. Zupan: {\it Classification of trisections and the Generalized Property R Conjecture},  Proc. Amer. Math. Soc. {\bf 144}(11), 4983-4997 (2016).

\bibitem{Miller-Naylor}
M. Miller, P. Naylor: {\it Trisections of non-orientable 4-manifolds}, Michigan Math. J.  {\bf 74}(2), 403-447 (2024). 

\bibitem{Montesinos}  J. M. Montesinos: {\it Heegaard diagrams for closed $4$-manifolds}, in:  {\em ``Geometric Topology"}, Academic Press, 1979. 
 

\bibitem{Rubinstein-Tillmann} 
J.~H. Rubinstein, S. Tillmann: {\it Multisections of piecewise linear manifolds},  Indiana Univ. Math. J.  {\bf 69}(6), 2209-2239 (2020). 

\bibitem{Spreer-Tillmann(Exp)} 
J. Spreer, S. Tillmann: {\it Determining the Trisection Genus of Orientable and Non-Orientable PL 4-Manifolds through Triangulations},  Exp. Math. {\bf  31}(3), 897-907 (2022).  

\bibitem{[Trace 1982]}
B. Trace: {\it On attaching 3-handles to a 1-connected 4-manifold}, Pacific Journal of Mathematics, {\bf 99}(1), 175-181 (1982).   

\end{thebibliography}
\end{document}